\documentclass[12pt]{article} %***
\usepackage[sectionbib]{natbib}
\usepackage{array,epsfig,fancyheadings,rotating}
\usepackage[]{hyperref}  %<----modified by Ivan
\usepackage{sectsty, secdot}
\sectionfont{\fontsize{12}{14pt plus.8pt minus .6pt}\selectfont}
\renewcommand{\theequation}{\thesection\arabic{equation}}
\subsectionfont{\fontsize{12}{14pt plus.8pt minus .6pt}\selectfont}
\usepackage[utf8]{inputenc}
\usepackage{amsmath}
\usepackage{amssymb}
\usepackage{amsfonts}
\usepackage{multirow}
\usepackage{amsthm}
\usepackage{natbib}
\usepackage{bm,amssymb,amsmath,amsthm,dsfont}
\usepackage{color}

\newtheorem{theorem}{Theorem}
\newtheorem{lemma}{Lemma}

\newtheorem{proposition}{Proposition}
\theoremstyle{definition}
\newtheorem{definition}{Definition}
\newtheorem{example}{Example}
\newtheorem{remark}{Remark}
\DeclareMathOperator{\diag}{diag}

\newcommand{\eeta}{\boldsymbol \eta}
\newcommand{\rrho}{\boldsymbol \rho}

\renewcommand{\H}[1]{\mathcal H_{#1}}

\newcommand{\aaaa}{\mathbf a}

\begin{document}

%%%%%%%%%%%%%%%%%%%%%%%%%%%%%%%%%%%%%%%%%%%%%%%%%%%%%%%%%%%%%%%%%%%%%%%%%%%%%%%%%%%%%%%%%%%%%%%%%%%%%%%%%%%%%%%%%%%%%%%%%%%%
%%%%%%%%%%%%%%%%%%%%%%%%%%%%%%%%%%%%%%%%%%%%%%%%%%%%%%%%%%%%%%%%%%%%%%%%%%%%%%%%%%%%%%%%%%%%%%%%%%%%%%%%%%%%%%%%%%%%%%%%%%%%

\renewcommand{\baselinestretch}{2}

\markright{ \hbox{\footnotesize\rm Statistica Sinica
%{\footnotesize\bf 24} (201?), 000-000
}\hfill\\[-13pt]
\hbox{\footnotesize\rm
%\href{http://dx.doi.org/10.5705/ss.20??.???}{doi:http://dx.doi.org/10.5705/ss.20??.???}
}\hfill }

\markboth{\hfill{\footnotesize\rm } \hfill}
{\hfill {\footnotesize\rm Identifiability of Bifactor Models} \hfill}

\renewcommand{\thefootnote}{}
$\ $\par

%%%%%%%%%%%%%%%%%%%%%%%%%%%%%%%%%%%%%%%%%%%%%%%%%%%%%%%%%%%%%%%%%%%%%%%%%%%%%%%%%%%%%%%%%%%%%%%%%%%%%%%%%%%%%%%%%%%%%%%%%%%%

\fontsize{12}{14pt plus.8pt minus .6pt}\selectfont \vspace{0.8pc}
\centerline{\large\bf Identifiability of Bifactor Models}
%\vspace{2pt} 
%\centerline{\large\bf}
\vspace{.4cm} 
\centerline{Guanhua Fang$^1$, Jinxin Guo$^2$, Xin Xu$^2$, Zhiliang Ying$^3$ and Susu Zhang$^4$} 
\vspace{.4cm} 
\centerline{\it $^1$ Baidu Research ~ 
	$^2$ Northeast Normal University}
\centerline{\it
$^3$ Columbia University, \\
$^4$ University of Illinois at Urbana-Champaign}
 \vspace{.55cm} \fontsize{9}{11.5pt plus.8pt minus.6pt}\selectfont

%%%%%%%%%%%%%%%%%%%%%%%%%%%%%%%%%%%%%%%%%%%%%%%%%%%%%%%%%%%%%%%%%%%%%%%%%%%%%%%%%%%%%%%%%%%%%%%%%%%%%%%%%%%%%%%%%%%%%%%%%%%%

\begin{quotation}
\noindent {\it Abstract:}
%{\bf Contents of the Abstract.}\\
The bifactor model and its extensions are multidimensional latent variable models, under which each item measures up to one subdimension on top of the primary dimension(s). Despite their wide applications to educational and psychological assessments, this type of multidimensional latent variable models may suffer from non-identifiability, which can further lead to inconsistent parameter estimation and invalid inference. The current work provides a relatively complete characterization of identifiability for the linear and dichotomous bifactor models and the linear extended bifactor model with correlated subdimensions. In addition, similar results for the two-tier models are also developed. Illustrative examples are provided on checking model identifiability through inspecting the factor loading structure. Simulation studies are reported that examine estimation consistency when the identifiability conditions are/are not satisfied.

\vspace{9pt}
\par
\end{quotation}\par

\def\thefigure{\arabic{figure}}
\def\thetable{\arabic{table}}

\renewcommand{\theequation}{\thesection.\arabic{equation}}

\fontsize{12}{14pt plus.8pt minus .6pt}\selectfont

\section{Introduction}

The bifactor method \citep{holzinger1937bi} for factor analysis is a constrained factor analytic model, which assumes that the responses to a set of test items can be accounted for by $(G+1)$ uncorrelated latent dimensions, with one primary dimension assessed throughout the test and $G$ secondary ``group'' dimensions. It further constrains each item to have nonzero loading on only one of the $G$ secondary dimensions. While the bifactor method for factor analysis (henceforth referred to as the linear bifactor model) was originally developed for continuous indicators, it has been extended to bifactor item response models \citep[e.g.,][]{gibbons1992full,gibbons2007full, cai2011generalized} for dichotomous, ordinal, or nominal item responses through the introduction of link functions, such as a probit link. The assumption of orthogonality among the secondary dimensions has also been relaxed in the extended bifactor model \citep[e.g.,][]{jennrich2012exploratory,jeon2013modeling} to allow for covariance between secondary dimensions unexplained by the primary dimension. The bifactor model with one primary dimension has further been extended to the two-tier model \citep{cai2010two} with $L\geq 1$ primary dimensions and $G$ secondary dimensions, with each item measuring up to one secondary dimension and the secondary dimensions being independent of the primary ones. 

The bifactor model and its extensions have demonstrated significant practical merits in educational and psychological assessments. Compared to uni- or low-dimensional latent trait models, they can accommodate the local dependence among clusters of items measuring the same subdimensions and produce subdimension trait estimates. Compared to general multidimensional latent variable models, they not only allow for the production of overall score(s) but also remarkably reduce the computational burden of high dimensional latent trait model estimation. The bifactor model and its extensions have hence been applied to hundreds of cognitive and psychological assessments, including psychiatric screenings that cover various domains of clinical disorders \citep[e.g.,][]{gibbons2009psychometric}, personality instruments that tap on multiple facets of the same trait \citep[e.g.,][]{chen2012modeling}, intelligence batteries with multiple subscales \citep[e.g.,][]{gignac2013bifactor}, and patient-reported outcome measures with broad situational representations \citep[e.g.,][]{reise2007role}. In educational testing, the bifactor model and its variants have seen wide applications to assessments that involve testlets, that is, multiple questions originated from the same stem \citep[e.g., passage; ][]{bradlow1999bayesian,demars2006application,demars2012confirming,jeon2013modeling,rijmen2010formal}. In longitudinal assessments with repeated administrations of the same item, bifactor and two-tier models can account for the within-person dependence of responses to the same item across time points \citep[see,][]{cai2016item}. Parameter estimation for the bifactor and two-tier models have been implemented in many commercial and open-source statistical software programs. 
{ It is also worthwhile to note that bifactor model is robust in practice. It tends to fit any data set better than other confirmatory models regardless of the population’s true models \citep{caspi2014p}.}
Thorough introductions to the bifactor model and its generalizations can be found in \cite{reise2012rediscovery} and \cite{cai2011generalized}.

Identifiability is a key issue in any type of latent variable modelling
\citep[][]{allman2009identifiability, xu2016identifiability, xu2017identifiability, gu2019sufficient, chen2015statistical, chen2019structured}. Intuitively, a model is identifiable if distinct parameter values produce unique probability distributions of observed responses. Model identifiability is a necessary condition for consistency of parameter estimation and valid statistical inference. Without any additional requirements, the bifactor model has been shown to be non-identifiable in previous studies: Under the linear bifactor model, \cite{green2018empirical} showed that two distinct sets of model parameters could produce the same model-implied covariance matrix. \cite{eid2018bifactor} showed that non-identifiability could arise in structural equation models with bifactor measurement model. 
The current paper addresses the identifiability issue of the bifactor model and its extensions by providing a relatively complete theory. We obtain the sufficient and necessary conditions for the identifiability of the standard bifactor model with both continuous and binary responses. Further, the necessary conditions for the extended bifactor model identifiability are given, and sufficient conditions for the identifiability of the extended bifactor and two-tier models are provided. For the dichotomous responses, discussions are limited to probit item response models, as the theoretical identifiability of logistic item response models with multivariate-normally distributed latent traits is a more complex issue. The identifiability of the aforementioned models can be achieved through the constraints on the loading structure of the items onto the general and specific dimensions. These results provide practitioners with viable means to examine identifiability of a certain test through a set of easily checkable conditions.

Many easily checkable conditions for general confirmatory factor model identifiability have been proposed in prior studies. Among them, the $t-$rule \cite[see][]{bollen1989structural}, which requires the number of unknown model parameters to not exceed the number of unique covariance terms, provides the necessary but not sufficient conditions for model identification. Another set of well-known rules are the 3-indicator and 2-indicator rules \cite[][]{bollen1989structural}, which are sufficient for identification but require simple factor loading structure, thus not applicable to the bifactor models. Empirical tests for local identification based on the information or the Jacobian matrix have been implemented in factor model estimation programs such as LISREL \cite[][]{joreskog2006lisrel}. However, local identifiability of parameters within the neighborhood of the estimates does not guarantee global identifiability in the entire parameter space. 
{ To the authors' knowledge, easily checkable sufficient and necessary conditions for general factor model identifiability are not yet available. }
Thanks to the special structure of the bifactor model, simple and checkable sufficient and necessary conditions can be developed. Readers are referred to \cite{bollen1989structural} and the LISREL manuals \cite{joreskog2006lisrel} for introductions to common methods for checking general factor model identification.

A related but distinct issue from model identification is rotational indeterminacy of exploratory factor models. Although methods to resolve rotational indeterminacy have been proposed \cite[][]{joreskog1979general}, it should be noted that such conditions do not guarantee confirmatory model identifiability \cite[][]{bollen1985uniqueness}. Although exploratory methods for bifactor analysis have been long-existing \cite[][]{holzinger1937bi, jennrich2012exploratory}, for the current study, discussions are restricted to confirmatory bifactor models where the loading structure is specified a priori, which are free from the issue of rotational indeterminacy.
%Identifiability conditions have been studied in previous research for both general confirmatory factor analytic models with continuous and ordinal responses and their special cases, such as the congeneric and multitrait-multimethod (MTMM) models \cite[e.g.,][]{reilly1995necessary,reilly1996identification,bollen2009two,shapiro1985identifiability,millsap2004assessing,grayson1994identification}. As another special case of the general factor analytic model with simple secondary loading structure and uncorrelated secondary dimensions, the bifactor model and its extensions can require weaker and more specific identifiability conditions than the general CFA model. Yet, to the authors' knowledge, existing research on the linear and dichotomous bifactor model and extended bifactor model has not discussed the minimal requirements for parameter identification, nor have there been established results on the contraints that can guarantee the identifiability of the two-tier model. 

The rest of the paper is organized as follows. Sections 2 presents the identifiability results on the linear bifactor, extended bifactor, and two-tier models. Section 3 extends the theoretical results on the three models to dichotomous responses with the probit link. Section 4 discusses the connections between the new results and the existing literature on bifactor identifiability.  A discussion of the findings is provided in Section 5. 
In the supplementary file, simulation studies are designed to verify the theoretical identifiability results by examining the estimation consistency under several identifiable and non-identifiable loading structures.
Additionally, the supplementary contains proofs of the main theoretical results and more detailed examples.

\section{Linear Bifactor Model and Extensions}

The current section presents results on the identifiability of the linear bifactor model and its extensions, which assume the response to each item to be normally distributed, with mean equal to an intercept plus a linear combination of the latent factor scores. This class of models is thus suitable for continuous observed indicators. 
%The first subsection discusses the results on the standard bifactor model, the second subsection pertains the extended bifactor model allowing correlation among subdimensions, and the third subsection provides the results on the two-tier model. For each model, the model setup and the definition of identifiability are introduced, the theoretical results on model identifiability are given, and examples are presented to illustrate how the theorems can be applied to determine parameter identifiability based on the loading matrices.

It is worth introducing the concept of identifiability in general before moving on to specific models. In mathematical terms, a statistical model may be specified by a pair $(\mathcal S, \mathcal P)$, where $\mathcal S$ is the set of possible observations or the sample space, and $\mathcal P$ is a set of probability distributions on $\mathcal S$, which is parameterized as $\mathcal P = \{P_{\theta}, \theta \in \Theta\}$. The set $\Theta$ defines the parameter space. We say the model parameter $\theta^{\ast}$ is identifiable (or model is identifiable at $\theta^{\ast}$) if $\mathcal P_{\theta}(y) = \mathcal P_{\theta^{\ast}}(y)$ for all $y \in \mathcal S$ implies $\theta = \theta^{\ast}$. Essentially, identifiability implies that the underlying distribution of the observed data cannot admit two distinct sets of parameter values.

\subsection{Standard bifactor model}\label{lin:bifactor}

The standard bifactor model \citep{holzinger1937bi} assumes that the response to each item in a test can be explained by one general factor, which runs through the test, and up to one group factor, which runs through a subset of items. Without loss of generality, we refer to a subset of items which load on the same group factor as a testlet \citep[e.g.,][]{bradlow1999bayesian, demars2006application}, but the model formulation is equally applicable to psychological assessments with subdimensions and cognitive batteries with subtests. Specifically, consider a test with $J$ items which can be partitioned into $G$ testlets, where the $g$th testlet consists of $J_g$ items and $\sum_g J_g = J$. Let $\mathcal{B}_g$ denote the set of items in the $g$th testlet. Under the standard linear bifactor model, the response to item $j \in \{1,\ldots, J\}$ in testlet $g_j \in \{1,\ldots, G\}$, $Y_j$, is given by 
\begin{eqnarray}\label{gaussmodelfull}
Y_j = d_j + a_{j0} \eta_0 + \sum_{g=1}^G a_{j g} \eta_{g} + \epsilon_{j},
\end{eqnarray}
where $\eta_0$ is the respondent's latent score on the general factor, $\eta_g$ is the latent score on the $g$th group factor, $a_{jk}$ is item $j$'s loading on the $k$th latent dimension, $d_j$ is the item intercept, and $\epsilon_{j}$ is the random error unexplained by the latent factors. Across all items, the $\epsilon_j$s are assumed to be independently distributed with mean $0$ and variance $\lambda_j$, that is, $\epsilon_j \sim N(0, \lambda_j)$. Further, for item $j$ in testlet $g_j$, it is assumed that $a_{jg} = 0$ for all $g\neq g_j$, in other words, the loadings of item $j$ on all other group factors are restricted to $0$. Thus, Equation \eqref{gaussmodelfull} simplifies to 
\begin{eqnarray}\label{gaussmodel}
Y_j = d_j + a_{j0} \eta_0 + a_{j g_j} \eta_{g_j} + \epsilon_{j}.
\end{eqnarray}
Let $\boldsymbol\eta = (\eta_0, \eta_1, \ldots, \eta_G)^T$ denote the vector of latent traits of a respondent, which is assumed to follow the multivariate normal distribution with zero mean and covariance $\boldsymbol{\Sigma}$, that is 
\begin{equation}
\boldsymbol{\eta} \sim \mbox{MVN}(\boldsymbol{0}, \boldsymbol{\Sigma}).
\end{equation}
The standard bifactor model further assumes that all general and group factors are independent, in other words, $\boldsymbol\Sigma = \boldsymbol{I}_{(1+G), (I+G)}$, where $\boldsymbol{I}_{(1+G), (I+G)}$ is the $(1+G)\times (1+G)$ identity matrix. The mean and standard deviation of each latent factor are fixed to $0$ and $1$, respectively, to resolve the location and scale indeterminacy of the latent dimensions.

Let $\boldsymbol A = [\mathbf a_0, \mathbf a_1, \ldots, \mathbf a_G]$ denote the $J\times (1+G)$ matrix of factor loadings, where the first column $\mathbf a_0 = (a_{10}, \ldots, a_{J0})^T$ is the items' loadings on the general factor, and the subsequent $G$ columns ($\mathbf a_1, \ldots, \mathbf a_G$) are the loadings on each of the $G$ testlets. Note that $\boldsymbol{A}$ is a sparse matrix, with most of the testlet-specific loadings restricted to zero by the single group factor loading assumption. Let $\mathbf d = (d_1, \ldots, d_J)^T$ denote the length-$J$ vector of item intercepts, and let $\boldsymbol\lambda = (\lambda_1, \ldots, \lambda_J)$ denote the vector of item unique variances. To resolve the sign indeterminacy of the latent factors, we further assume without loss of generality that for each factor, the first item that loads on the factor has positive loading. Under the standard bifactor model, the parameter space $\mathfrak P$ is then given by $\mathfrak P = \{(A, \mathbf d, \boldsymbol \lambda) ~ |~ \textrm{The sign of first non-zero element in every column $ $ of $A$ is positive} \}.$ The definition of identifiability of the standard bifactor model is as follows.

\begin{definition}
	We say a linear bifactor model is identifiable at $(A, \mathbf d, \boldsymbol \lambda)$ if, for any other set of parameter $(A^{'}, \mathbf d^{'}, \boldsymbol \lambda^{'})$ that defines the same probability distribution, it must hold that 
	\begin{eqnarray}
	(A, \mathbf d,  \boldsymbol \lambda) = (A^{'}, \mathbf d^{'}, \boldsymbol \lambda^{'}).
	\end{eqnarray} 
\end{definition}

Before stating the identifiability results on the standard bifactor model, we introduce some additional notation. 
%Let $A[l,\cdot]$ and $A[\cdot,l]$ denote the $l$th row and $l$th column of matrix $A$, respectively, and let $A[-l,\cdot]$ (or $A[\cdot,-l]$) denote the submatrix of $A$ by removing the $l$th row (or $l$th column). Similar, for any subsets $\mathcal{B}_1,\mathcal{B}_2$, let $A[\mathcal{B}_1, \mathcal{B}_2]$ denote the submatrix of $A$ consisting only of rows in $\mathcal{B}_1$ and columns in $\mathcal{B}_2$, and let $A[-\mathcal{B}_1,\cdot]$ be the submatrix of $A$ removing rows in $\mathcal B_1$ (similarly for the columns). 
Let $\bar A_g$ be the submatrix of $A$ corresponding to items in testlet $g$ (i.e., $B_g$), that is, $ \bar A_g = (\mathbf a_{0}[\mathcal B_g], \mathbf a_{g}[\mathcal B_g])$, where $x[\mathcal{B}_g]$ denotes of subvector with entries in set $\mathcal{B}_g$. Further, let $A[\mathcal{B}, :]$ denote the submatrix of $A$ consisting only of rows in some generic set $\mathcal{B}$. We define the following subsets. 
\begin{itemize}
	\item $\H{1} = \{g ~| ~ \mathbf a_{0}[\mathcal B_g] \neq \mathbf 0 \}$: the set of testlets in \{1,\ldots, G\} with non-vanishing main factors, that is, at least one item in the testlet has nonzero true loading on the main dimension, $\eta_0$;
	% \item $\H{3} = \{g ~| ~ \bar A_g ~\textrm{has column rank 2} \}$: the set of testlets with full main factor and testlet-specific factor information, that is, columns $\mathbf a_{0}[\mathcal B_g]$  and $\mathbf a_{g}[\mathcal B_g]$ are linearly independent;
	% \item $\mathcal N = \{g ~| ~ \mathbf a_g[\mathcal B_g] = \mathbf 0 \}$: the set of testlets with vanishing testlet-specific factors, that is, items in the testlet do not load on any group dimension.
	\item $\mathcal Q_{g} = \{j ~ | ~ \mathbf a_g[j] \neq 0 \}$: the set of items in the $g$th testlet with nonzero true loadings on the testlet-specific factor;
	\item $\H{2} = \{g ~ | ~\textrm{there exists a partition of } \mathcal B_g, i.e., \mathcal B_g = \mathcal B_{g,1} \cup \mathcal B_{g,2},~ \mathcal B_{g,1} \cap \mathcal B_{g,2} = \emptyset, \textrm{ such that } \bar A_g[\mathcal B_{g,1},:], \bar A_g[\mathcal B_{g,2},:]~\textrm{are of full column rank.} \}$: the set of testlets that can be partitioned into two disjoint subsets of items, where each submatrix of (main and testlet) factor loadings has full column rank.
\end{itemize} 

The following theorem characterizes the sufficient and necessary conditions for the identifiability of the linear standard bifactor model. 
\begin{theorem}\label{thm:unknown:bifactor}
	Under the standard linear bifactor model, the model parameters are identifiable if and only if it satisfies one of the following conditions:
	\begin{itemize}
		\item[P1~~] $|\mathcal Q_g| \geq 3$ for all $g = 1, \ldots, G$; $|\H{1}| \geq 3$.
		\item[P2~~] $|\mathcal Q_g| \geq 3$ for all $g = 1, \ldots, G$; $|\H{1}| = 2 $; $|\H{2}| \geq 1$.
	\end{itemize}
\end{theorem}

Theorem \ref{thm:unknown:bifactor} gives the minimum requirements for identifiability of the standard bifactor model. Specifically, it requires the test to contain at least $2$ testlets, each containing at least $3$ items. In addition, if there are only $2$ testlets, Theorem \ref{thm:unknown:bifactor} requires that one of them can be partitioned into $2$ disjoint subsets of items, such that both subsets have linearly independent primary and testlet-specific factor loadings.

%\begin{example}
%	Consider a standard bifactor model for $J = 7$ items, where the true parameters are
%	\begin{eqnarray}
%	A = 
%	\begin{pmatrix}
%	a & b & 0 & 0 \\
%	a & b & 0 & 0 \\
%	c & d & 0 & 0 \\
%	c & d & 0 & 0 \\
%	e & 0 & f & 0 \\
%	e & 0 & f & 0 \\
%	e & 0 & f & 0 \\
%	\end{pmatrix},
%	\quad \mathbf d = 
%	\begin{pmatrix}
%	d_1 \\
%	d_2 \\
%	d_3 \\
%	d_4 \\
%	d_5 \\
%	d_6 \\
%	d_7 \\
%	\end{pmatrix},
%	\end{eqnarray}
%	with $a,b,c,d,e,f \neq 0$ and $a d \neq b c$. One can check that (1) $|\H{1}| = 2$, because both testlet $1$ and testlet $2$ have nonzero main factor loadings; (2) $|\mathcal Q_g|\geq 3$ for all $g$, because $|\mathcal Q_1| = 4, |\mathcal Q_2| = 3$; and (3) $|\H{2}|\geq 1$, because testlet $1$ can be partitioned to $\mathcal B_{g,1} = \{1,3\}, \mathcal B_{g,2} = \{2,4\}$, each containing linearly independent columns. By Condition P2, we know that the model is identifiable.
%\end{example}

\subsection{Extended bifactor model}\label{lin:extended}
The extended bifactor model, also known as the oblique bifactor model \citep{jennrich2012exploratory}, relaxes the assumption of independence between the secondary dimensions. Instead of restricting the latent covariance matrix, $\boldsymbol{\Sigma}$, to be the identity matrix, the extended bifactor model allows the covariance between the latent dimensions to take the form of $\boldsymbol\Sigma = \begin{pmatrix}
1& \mathbf 0^T \\
\mathbf 0^T  & \boldsymbol\Sigma_G \\
\end{pmatrix}$, where the covariance matrix for the testlet dimensions, $\boldsymbol\Sigma_G$, is positive-definite with all diagonal elements being 1 and no additional restriction on off diagonal elements. Note that the covariances between the primary dimension and each testlet dimension is still restricted to be $0$. 

Under the extended bifactor model, in addition to the item intercepts, loadings, and unique variances, the latent covariance matrix also needs to be estimated.  
Let $\mathfrak P = \{(A, \mathbf d, \Sigma_G, \boldsymbol \lambda) ~ |~  \textrm{first non-zero element in every column}$ $ \textrm{of $A$ is positive}, \mathrm{diag}(\Sigma_G) = 1, \Sigma_G ~ \textrm{positive definite} \}$ denote the parameter space of the extended bifactor model. Then the identifiability of the extended bifactor model is defined as follows.

\begin{definition}
	We say a linear extended bifactor model is identifiable at $(A, \mathbf d, \Sigma_G, \boldsymbol \lambda)$ if for any other set of parameters $(A^{'}, \mathbf d^{'}, \Sigma_G^{'}, \boldsymbol \lambda^{'})$ that defines the same probability distribution, 
	\begin{eqnarray}
	(A, \mathbf d, \Sigma_G, \boldsymbol \lambda) = (A^{'}, \mathbf d^{'}, \Sigma_G^{'}, \boldsymbol \lambda^{'})
	\end{eqnarray} 
	must hold. 
\end{definition}

In addition to $\H{2},$ and $\mathcal Q_g$ in the standard bifactor model identifiability results, we introduce another set that is key to the identifiability of the extended bifactor model:
\begin{itemize}
	\item $\H{3} = \{g ~| ~ \bar A_g ~\textrm{has column rank 2} \}$: the set of testlets with full main factor and testlet-specific factor information, that is, columns $\mathbf a_{0}[\mathcal B_g]$  and $\mathbf a_{g}[\mathcal B_g]$ are linearly independent.
\end{itemize}

Theorem \ref{thm:extendunknown} gives two sets of sufficient conditions for identifiability of the extended linear bifactor model.
\begin{theorem}\label{thm:extendunknown}
	Under the linear extended bifactor model, the model parameters are identifiable if one of the following sets of requirements is satisfied:
	\begin{itemize}
		\item [E1S~] $|\mathcal Q_g| \geq 3$ for all $g = 1, \ldots, G$; $|\H{3}| \geq 3$.
		\item [E2S~] $|\mathcal Q_g| \geq 3$ for all $g = 1, \ldots, G$; $|\H{3}| = 2$; $|\H{2}| \geq 1$. 
	\end{itemize}
\end{theorem}

The sufficient conditions in Theorem \ref{thm:extendunknown} are very similar to the sufficient and necessary conditions for the standard bifactor model in Theorem \ref{thm:unknown:bifactor}, where $\H{3}$ is the counterpart to $\H{1}$. Unlike $\H{1}$, which contains testlets with nonzero main-factor loading vector, $\H{3}$ further requires the main- and testlet-factor loading vectors to be linearly independent. 
Note that $\mathit{E1S}$ and $\mathit{E2S}$ are sufficient for identifiability of the extended bifactor model but they are not necessary. Theorem \ref{thm:extendunknown:nec} provides the necessary but not sufficient conditions, that is, the minimum conditions that need to be met.

\begin{theorem}\label{thm:extendunknown:nec}
	Under the linear extended bifactor model, the model parameters are identifiable only if both conditions below are satisfied.
	\begin{itemize}
		\item [E1N~] $|\mathcal Q_g| \geq 2$ for all $g = 1, \ldots, G$; and $|\mathcal Q_g| \geq 3$ for $g$ such that $\Sigma_G[g,-g] = \mathbf 0$. 
		\item [E2N~] $|\H{3}| \geq 2$. 
	\end{itemize}
\end{theorem}

Essentially, at least two testlets with linearly independent main- and testlet-factor loadings are required. In addition, each subdimension should be measured by at least $2$ items, and if a subdimension is uncorrelated with others, at least $3$ items are required as before. The nonzero correlation between the testlet factors provides additional information on the testlet-specific loadings, reducing the number of required items per testlet to $2$ for those testlets with nonzero correlations with others. 

Note that $\mathit{E1N}$ and $\mathit{E2N}$ together are not enough for the identifiability of extended bifactor model. Additional requirements are needed to render parameter identifiability. 
% Do not remove. The construction of counter example when only E1N and E2N hold. 
%\begin{eqnarray}
%\sigma^{'} = \mathrm{sgn}(\sigma) \sqrt{\frac{ac \sin^2 \theta - \sigma \cos \theta \sin \theta (bc + ad) + bd \sigma^2 \cos^2 \theta}{ac \sin^2 \theta + bd (1 - \sigma^2 \sin^2 \theta) + \sigma \cos \theta \sin \theta (bc + ad)}}.
%\end{eqnarray}
Theorem \ref{thm:extendunknown} gives one way to impose such additional requirements. The requirement of $|\H{2}| \geq 1$ in $\mathit{E2S}$ may be replaced by other requirements. See the following proposition.

\begin{proposition}\label{prop:extend:H5:unnec}
	Under the linear extended bifactor model, the model parameters are identifiable if $E1N$ and $E2N$ are satisfied, and
	\begin{itemize}
		\item[E3S] There exists $g_1, g_2 \in \H{3}$, such that (1) $\Sigma_G[g_1,g_2] \neq 0$; (2) $|\mathcal Q_{g_1}| \geq 3, |\mathcal Q_{g_2}| \geq 3$; and (3) Kruskal rank of $\bar A_{g_1}^T$ is 2. (A matrix $A$ has Kruskal rank \citep{kruskal1977three} $R$ if any $R$ columns of $A$ are linearly independent).   
		%\sz{$\bar A_{g_1}^T$ or $\bar A_{g_1}$? Proof says $\bar A_{g_1}$ has Kruscal rank $2.$ Also, do both $g_1$ and $g_2$ need exactly $3$ items?}  
	\end{itemize}
\end{proposition}

\vspace{-5mm}

\begin{remark}
The gap between the necessary and sufficient conditions in Theorems \ref{thm:extendunknown} and \ref{thm:extendunknown:nec} is that $|\mathcal H_3| \geq 2$ itself cannot guarantee identifiability. Either more testlets ($|\mathcal H_3| \geq 3$) are needed or at least one of the testlets needs to be ``strong'' ($|\mathcal H_2| \geq 1$). Non-zero correlations ($\Sigma$) between testlet factors increase the complexity of the identifiability problem compared to the uncorrelated case.  
\end{remark}

\vspace{-5mm}

	\begin{remark}
Conditions in Theorems \ref{thm:unknown:bifactor} - \ref{thm:extendunknown:nec} are easy to check in practice. We only need to do simple algebras (i.e. counting the number of non-zero entries, computing the column rank, etc.) on the estimated loading matrix $A$ and covariance matrix $\Sigma$. 
Additionally, $\mathcal H_2 \geq 1$ in Condition P2 generically holds when a testlet contains four or more items. 
	\end{remark}

\vspace{-5mm}

\begin{example}
	Consider an extended bifactor model with three testlets, where testlet 1 has only two items. Suppose the true $A$ and $\Sigma_G$ are given by 
	%\sz{changed notations to be consistent with setup; Also, it seems that this example is not matched to either of the sufficiency theorems (Theorem \ref{thm:extendunknown} or Proposition \ref{prop:extend:H5:unnec})?}
	\begin{eqnarray}
	A = 
	\begin{pmatrix}
	a_{10} & a_{11} & 0      & 0 \\
	a_{20} & a_{21} & 0      & 0 \\
	a_{30} & 0      & a_{32} & 0 \\
	a_{40} & 0      & a_{42} & 0 \\
	a_{50} & 0      & a_{52} & 0 \\
	a_{60} & 0      & a_{62} & 0 \\
	a_{70} & 0      & 0      & a_{73} \\
	a_{80} & 0      & 0      & a_{83} \\
	a_{90} & 0      & 0      & a_{93} \\
	\end{pmatrix}, ~~~~~~
	\Sigma_G = 
	\begin{pmatrix}
	1 & \sigma_{12} & \sigma_{13} \\
	\sigma_{12} & 1 &\sigma_{23} \\
	\sigma_{13} & \sigma_{23} & 1 \\
	\end{pmatrix},
	\end{eqnarray}
	where $\sigma_{12}, \sigma_{13}, \sigma_{23} \neq 0$, any $A_g$ {($g = 1, 2, 3$)} has 2 linearly independent columns.
	According to Proposition \ref{prop:extend:H5:unnec}, the model parameter is identifiable, even though $|\mathcal Q_1| = 2$.
\end{example}

\subsection{Two-tier model}\label{lin:twotier}

The two-tier model \citep{cai2010two} extends the standard bifactor model by allowing for more than one primary dimensions. Consider a test that measures $L$ primary factors and $G$ group factors. Under the two-tier model, denote latent factors by $\eeta = (\eeta_1, \eeta_2)$ with $\boldsymbol \eta_1 = (\eta_1, \ldots, \eta_L)^T$ and $\boldsymbol \eta_2 = (\eta_{L+1}, \ldots, \eta_{L+G})^T$. The response to the $j$th item in testlet $g_j$ is given by
\begin{eqnarray}
Y_j = d_j + \sum_{l = 1}^L a_{jl} \eta_{l} + \sum_{g=L+1}^{L+G} a_{jg} \eta_{g} + \epsilon_j.
\end{eqnarray}
Similar to the bifactor model, $\epsilon_j$s are indepdent and normally distributed with $\epsilon_j\sim N(0, \lambda_j)$, and only $g_j$th testlet factor loading is nonzero for item $j$, i.e., $a_{jg}=0$ for all $g\neq g_j$. The latent covariance matrix of $\eeta$ takes the form of $\begin{pmatrix}
\Sigma_L& \mathbf 0^T \\
\mathbf 0^T  & I_{G \times G} \\
\end{pmatrix}$, where $\Sigma_L$ is a $L\times L$ positive-definite matrix with diagonal elements of 1 and no additional restriction on off-diagonal elements. Let  $\mathfrak P = \{(A, \mathbf d, \Sigma_L, \boldsymbol \lambda) ~ |~ \textrm{first non-zero element in every column of $A$ is positive}$,
$\mathrm{diag}(\Sigma_L) = 1, \Sigma_L ~ \textrm{positive definite} \}$ represent the model space of the two-tier model. Further, let  $\mathbf I$ denote the configuration mapping function, with $\mathbf I(X) = \tilde X$ where $\tilde X_{ij} = \mathbf 1\{X_{ij} \neq 0 \}$ for an arbitrary matrix $X$.
The definition of two-tier model identifiability is stated as follows.
\begin{definition}
	A linear two-tier model is identifiable at $(A, \mathbf d, \Sigma_L, \boldsymbol{\lambda})$ if for any other set of parameters $(A^{'}, \mathbf d^{'}, \Sigma_L^{'}, \boldsymbol{\lambda}^{'})$ that define the same distribution and satisfies $\mathbf I(A) = \mathbf I(A^{'})$, it must hold that $(A, \mathbf d, \Sigma_L, \boldsymbol{\lambda}) = (A^{'}, \mathbf d^{'}, \Sigma_L^{'}, \boldsymbol{\lambda}^{'})$. 
\end{definition}
Here, notice that a factor model is only identifiable up to some rotation. In the definition, the requirement for equal factor loading configurations, $\mathbf I(A^{'}) = \mathbf I(A)$, is put in place to resolve the rotational indeterminacy. We would like to point out that the identifiability of two-tier model is non-trivial, in the sense that the model could fail to be identifiable, even if (1) the loading matrix of main factors, $A_{:,1:L}$, satisfies the usual identifiability conditions for multivariate factor models and (2) the testlets satisfy the identifiability conditions for the bifactor model. See the examples provided in the supplementary.

The proposed sufficient conditions for two-tier model identifiability build upon the sufficient conditions for unique variance identifiability under factor models with uncorrelated errors, which can be found in Theorem 5.1 of \cite{anderson1956statistical} and is rephrased below. 
\begin{theorem}\citep{anderson1956statistical}.\label{thm:generalother}
	Consider a general factor model with implied covariance matrix $\Psi = A\Sigma A^T + \Lambda$, where the item error covariance matrix, $\Lambda$, is diagonal with $\diag(\Lambda) = \boldsymbol{\lambda}$. Then $A \Sigma A^T$ and $\Lambda$ are identifiable if the following holds:
	\begin{itemize}
		\item[C0] If any row of $A$ is deleted, there remain two disjoint submatrices of $A$ with full column rank.
	\end{itemize}
\end{theorem}

Under the two-tier model, let $\bar A_g$ be the submatrix of $A$ corresponding to items in testlet $g$, that is, $\bar A_g= (A[\mathcal B_g,1:L], A[\mathcal B_g, L+g])$. Similarly, for a subset of testlets $\mathcal G_1 \subseteq \{1,\ldots, G\}$, denote the submatrix of $A$ corresponding to testlets in $\mathcal G_1$ by $\bar A_{\mathcal G_1} = (A[\mathcal B_{\mathcal G_1},1:L],A[\mathcal B_{\mathcal G_1},L+\mathcal G_1])$, where $\mathcal B_{\mathcal G_1} = \bigcup_{g\in G_1}\mathcal B_g$ and $L+\mathcal G_1 = \{L+g^*\mid g^*\in \mathcal G_1\}$. We introduce two sets, $\H{4}$ and $\H{5}$, that are essential to the identifiability results for the two-tier model:
\begin{itemize}
	\item $\H{4} = \{g ~ | ~ A[\mathcal B_g, 1:L] ~\textrm{is of full column rank} \}$: the set of testlets with non-degenerate main factor information, that is, the $L$ columns corresponding to the main factor loadings are linearly independent;
	\item $\H{5} = \{g ~ | ~ \bar A_g ~\textrm{is of full column rank} \}$: the set of testlets with non-degenerate all-factor information, that is, with linearly independent main factor and testlet factor loadings.
\end{itemize}

The following theorem provides the sufficient conditions for identifiability of the linear two-tier model. 

\begin{theorem}\label{twotier:unknown}
	Under the linear two-tier model, if the true loading matrix, $A$, satisfies Condition C0 and one of Conditions T1S - T3S, then the parameters are identifiable.
	\begin{itemize}
		\item[T1S~] $|\H{4}| \geq 3$, $A[\mathcal B_{\H{4}}, 1:L]$ contains an identity, where $\mathcal B_{\H{4}}$ is the set of items that makes up the testlets in $\H{4}$.
		\item[T2S~] $|\H{4}| \geq 2$, $|\H{5}| \geq 1$, $A[\mathcal B_{\H{4}}, 1:L]$ contains an identity.
		\item[T3S~] A[:,1:L] contains an identity, and there exists a partition of testlets $ \{1, \ldots, G\} = \mathcal G_1 \dot\cup \mathcal G_2$, such that (1) $\bar A_{\mathcal G_1}$ has full-column rank, and
		(2) $A[\mathcal B_{\mathcal G_2},1:L]$ has full-column rank. 
	\end{itemize}
\end{theorem}

In other words, the linear two-tier model is identifiable if the following two conditions are simultaneously met:
\begin{itemize}
	\item Removing any row of the loading matrix $A$, the remaining rows of $A$ can be partitioned into two disjoint submatrices, both of which contain $L+G$ linearly independent columns.
	\item One of the following is satisfied:
	\begin{itemize}
		\item[$\mathit{T1S:}$] (1) The test contains at least $3$ testlets with linearly independent main factor loadings, and (2) within these testlets satisfying (1), for each main factor, there exists at least one item that exclusively measures this main factor (i.e., having nonzero loadings only on this main factor and possibly the testlet factor);
		\item[$\mathit{T2S:}$] (1) The test contains $2$ testlets with linearly independent main factor loadings, (2) within these testlets satisfying (1), for each main factor, there exists at least one item that exclusively measures this main factor, and (3) at least one of the testlets satisfying (1) has linearly independent main \textit{and} testlet factor loadings;
		\item[$\mathit{T3S:}$] (1) For each main factor, at least one item in the test exclusively measures that main factor (aside from the testlet factor), and (2) the set of all testlets 
		%($\{1, \ldots, G\}$) 
		can be partitioned to two disjoint subsets, $\mathcal G_1$ and $\mathcal G_2$, such that (a) the loading matrix corresponding to the first subset of testlets, $\bar A_{\mathcal G_1}$, has $L+|G_1|$ linearly independent columns, and (b) for the second subset, the columns of main factor loadings, $A[\mathcal B_{\mathcal G_2}, 1:L]$, are linearly independent.
	\end{itemize}
\end{itemize}

\section{Extensions to Dichotomous Responses}

Notwithstanding the wide application of linear factor models in the social science literature, a large proportion of educational and psychological assessments consist of items with dichotomous responses, for instance, cognitive questions where an examinee responds either correctly ($1$) or incorrectly ($0$), or clinical screening questions where a participant either exhibits certain behavior ($1$) or not ($0$). Probit item response models, also known as normal ogive models \citep[e.g.,][]{thurstone1927unit,lawley1943xxiii,lord1952theory,christoffersson1975factor}, have been widely adopted for dichotomous responses. In general, consider a test of $J$ items. The responses to the $J$ items are assumed to be locally independent given the respondents' latent traits, $\eeta$, and the probability of responding ``1'' to the $j$th item is given by 
\begin{equation}\label{probit}
P(Y_j = 1 | \eeta) = P_{\epsilon}(\epsilon_j \leq d_j +  \mathbf a_j^T \eeta) = \Phi(d_j + \mathbf a_j^T \eeta),
\end{equation}
where $\epsilon_j \sim N(0,1)$, $\Phi(\cdot)$ is the standard normal cumulative distribution function (i.e., the probit link), and $\eeta\sim N(0, \Sigma), d_j, \mathbf a_j$ are the person latent traits, item intercept, and item slopes/loadings, respectively. Similar as before, let $A$ denote the matrix of factor loadings.

Item bifactor model and extensions have also been proposed to accommodate for dichotomous response tests with underlying bifactor-like latent structures. Examples of such models include the item bifactor model \citep{gibbons1992full}, the extended item bifactor model \citep{jeon2013modeling}, and the two-tier item factor model \citep{cai2010two}. The identifiability results for the linear bifactor model and extensions do not directly apply to dichotomous item bifactor-type models due to different parameterizations and form of observed data. The current section presents the results on identifiability of the dichotomous bifactor model, extended bifactor model, and two-tier model with probit links. 

Before introducing the identifiability conditions for each of the specific models, it is worth mentioning a few identities on the items' $1$st, $2$nd, and $k$th moments implied by the general probit item factor model in \eqref{probit}, as well as their relationships with the thresholds and tetrachoric correlations \citep{pearson1900mathematical} under the probit model. Unlike linear models, we do not directly observe the mean and covariance matrix implied by the linear component (i.e., $d_j + \aaaa_j^T\eeta_j$) under the probit model, but the threshold and tetrachoric correlations can be identified \citep{kendall1958chaps} and estimated with various approximation methods \citep[e.g.,][]{castellan1966estimation, olsson1979maximum}. In the following, we explain how the tetrachoric correlations relate to the idenfiability problem.

Let $\xi$ denote a standard normal random variable. Note that, at the population level (i.e., for a randomly chosen $\eeta\sim MVN(0, \Sigma$)), the marginal probability of observing a response of $1$ on item $j$ is given by \begin{eqnarray}\label{1dim}
P(Y_j = 1) & = & \mathbb E_{\eeta} P(Y_j = 1 | \eeta) = 
\mathbb E_{\eeta} \mathbb E_{\epsilon_j} [\mathbf 1\{ \epsilon_j \leq d_j + \aaaa_j^T \eeta \} | \eeta]  =P(\epsilon_j \leq d_j + \aaaa_j^T \eeta) \nonumber \\
& = & P(d_j + \sqrt{\aaaa_j^T \Sigma \aaaa_j + 1} \xi \geq 0) = 1 - \Phi(- \frac{d_j}{\sqrt{\aaaa_j^T \Sigma \aaaa_j + 1}}).
\end{eqnarray}
And the probability that the responses to items $j_1$ and $j_2$ are both $1$ is 
\begin{eqnarray}\label{2dim}
& & P(Y_{j_1} = 1, Y_{j_2} = 1) 
=  \mathbb E_{\eeta} P(Y_{j_1} = 1, Y_{j_2} = 1 | \eeta)  \nonumber \\
% = \mathbb E_{\eeta}[P(Y_{j_1} = 1 | \eeta) P(Y_{j_2} = 1 | \eeta)]
% & = & \mathbb E_{\eeta} \{\mathbb E_{\epsilon_{j_1}} [\mathbf 1\{ \epsilon_j \leq d_j + a_j^T \eeta \} | \eeta] \mathbb E_{\epsilon_{j_2}} [\mathbf 1\{ \epsilon_j \leq d_j + a_j^T \eeta \} | \eeta] \} \nonumber \\
& = & P(\epsilon_{j_1} \leq d_{j_1} + \aaaa_{j_1}^T \eeta, \epsilon_{j_2} \leq d_{j_2} + \aaaa_{j_2}^T \eeta) \\
& = & P(d_{j_1} + \sqrt{\aaaa_{j_1}^T \Sigma \aaaa_{j_1} + 1} \xi_{j_1} \geq 0 , d_{j_2} + \sqrt{\aaaa_{j_2}^T \Sigma \aaaa_{j_2} + 1} \xi_{j_2} \geq 0) \nonumber \\
& = & \Phi_2 (- \frac{d_{j_1}}{\sqrt{\aaaa_{j_1}^T \Sigma \aaaa_{j_1} + 1}}, - \frac{d_{j_2}}{\sqrt{\aaaa_{j_2}^T \Sigma \aaaa_{j_2} + 1}}, \frac{\aaaa_{j_1}^T \Sigma \aaaa_{j_2}}{\sqrt{\aaaa_{j_1}^T \Sigma \aaaa_{j_1} + 1} \sqrt{\aaaa_{j_2}^T \Sigma \aaaa_{j_2} + 1}} ), \nonumber
\end{eqnarray}
where $\Phi_2(a,b, \rho) = \mathbb E(X_1 \geq a, X_2 \geq b)$, $X_1, X_2 \sim N(0,1)$ and $\mathrm{corr}(X_1, X_2) = \rho$. The $a,b$ and $\rho$ are commonly referred to as the thresholds and tetrachoric correlation under the probit framework.
The probability for responding $1$ simultaneously on $k$ items ($j_1, \ldots, j_k$) is given by
\begin{eqnarray}\label{kdim}
& & P(Y_{j_1} = 1, \ldots, Y_{j_k} = 1) 
 =  \mathbb E_{\eeta} P(Y_{j_1} = 1, \ldots Y_{j_k} = 1 | \eeta)  \nonumber \\
% = \mathbb E_{\eeta}[P(Y_{j_1} = 1 | \eeta) P(Y_{j_2} = 1 | \eeta)]
% & = & \mathbb E_{\eeta} \{\mathbb E_{\epsilon_{j_1}} [\mathbf 1\{ \epsilon_j \leq d_j + a_j^T \eeta \} | \eeta] \mathbb E_{\epsilon_{j_2}} [\mathbf 1\{ \epsilon_j \leq d_j + a_j^T \eeta \} | \eeta] \} \nonumber \\
& = & P(\epsilon_{j_1} \leq d_{j_1} + \aaaa_{j_1}^T \eeta, \ldots, \epsilon_{j_k} \leq d_{j_k} + \aaaa_{j_k}^T \eeta) \nonumber \\
& = & P(d_{j_1} + \sqrt{\aaaa_{j_1}^T \Sigma \aaaa_{j_1} + 1} \xi_{j_1} \geq 0 , \ldots,  d_{j_k} + \sqrt{\aaaa_{j_k}^T \Sigma \aaaa_{j_k} + 1} \xi_{j_k} \geq 0) \nonumber \\
& = & \Phi_k (- \frac{d_{j_1}}{\sqrt{\aaaa_{j_1}^T \Sigma \aaaa_{j_1} + 1}}, \ldots, - \frac{d_{j_2}}{\sqrt{\aaaa_{j_2}^T \Sigma \aaaa_{j_2} + 1}}, C_{\rho})
\end{eqnarray}
with tetrachoric correlation matrix $C_{\rho}[j_1,j_2] = \frac{\aaaa_{j_1}^T \Sigma \aaaa_{j_2}}{\sqrt{\aaaa_{j_1}^T \Sigma \aaaa_{j_1} + 1} \sqrt{\aaaa_{j_2}^T \Sigma \aaaa_{j_2} + 1}}$ and $C_{\rho}[j,j] = 1$.
Here $\Phi_k(a_1, \ldots ,a_k, C_\rho) = \mathbb E(X_1 \geq a_1, X_k \geq a_k)$, $X_1, \ldots, X_k \sim N(0,1)$ and $\mathrm{corr}(X_{k_1}, X_{k_2}) = C_\rho[k_1,k_2]$.

In the following, we show that threshold and tetrachoric correlations provide full information on probit binary item responses.
\begin{proposition}\label{prop:tetra}
	Two sets of parameters define the same model if and only if their thresholds and tetrachoric correlations are equal, i.e.,
	\begin{eqnarray}
	\frac{d_j}{\sqrt{\aaaa_j^T \Sigma \aaaa_j + 1}} = \frac{d_j^{'}}{\sqrt{(\aaaa_j^{'})^T \Sigma^{'} \aaaa_j^{'} + 1}} \quad \forall j, \nonumber 
	\end{eqnarray}  
	and
	\begin{eqnarray}
	\frac{\aaaa_{j_1}^T \Sigma \aaaa_{j_2}}{\sqrt{\aaaa_{j_1}^T \Sigma \aaaa_{j_1} + 1} \sqrt{\aaaa_{j_2}^T \Sigma \aaaa_{j_2} + 1}} = 
	\frac{(\aaaa_{j_1}^{'})^T \Sigma^{'} \aaaa_{j_2}^{'}}{\sqrt{(\aaaa_{j_1}^{'})^T \Sigma^{'} \aaaa_{j_1}^{'} + 1} \sqrt{(\aaaa_{j_2}^{'})^T \Sigma^{'} \aaaa_{j_2}^{'} + 1}} \quad \forall j_1 \neq j_2. \nonumber 
	\end{eqnarray} 
\end{proposition}

%\begin{proof}
%	The sufficient part is straightforward by noticing that $P(Y_{j_1} = 1, \ldots, Y_{j_k} = 1)$ only depends on $d_j / (\aaaa_j^T \Sigma \aaaa_j + 1)^{1/2}$'s and $ (\aaaa_{j_1}^T \Sigma \aaaa_{j_2}) /\big((\aaaa_{j_1}^T \Sigma \aaaa_{j_1} + 1)(\aaaa_{j_2}^T \Sigma \aaaa_{j_2} + 1)\big)^{1/2}$'s for all possible combinations of $j_1, \ldots  j_k$.
%	
%	The necessary part is also not hard. Notice that CDF function $\Phi$ is a strictly monotone increasing function. By \eqref{1dim}, we must have  
%	$d_j/(\aaaa_j^T \Sigma \aaaa_j + 1)^{1/2} = d_j^{'}/((\aaaa_j^{'})^T \Sigma^{'} \aaaa_j^{'} + 1)^{1/2}$ for all $j$.
%	In addition, $\Phi_2(a,b,\rho)$ is a strictly monotone increasing function of $\rho$ for any fixed $a, b$. Thus, from \eqref{2dim}, we get 
%	$(\aaaa_{j_1}^T \Sigma \aaaa_{j_2})/\big((\aaaa_{j_1}^T \Sigma \aaaa_{j_1} + 1) (\aaaa_{j_2}^T \Sigma \aaaa_{j_2} + 1)\big)^{1/2} = 
%	(\aaaa_{j_1}^T \Sigma^{'} \aaaa_{j_2})/\big(((\aaaa_{j_1}^{'})^T \Sigma^{'} \aaaa_{j_1}^{'} + 1) ((\aaaa_{j_2}^{'})^T \Sigma^{'} \aaaa_{j_2}^{'} + 1)\big)^{1/2}$ for any $j_1 \neq j_2$. Hence we prove the claim.
%\end{proof}

It follows from the above proposition that checking the identifiability of probit bifactor models comes down to checking whether the probit threshold and tetrachoric correlations admit only one set of parameters. In other words, the probit bifactor models can be identified if $(d_j, \aaaa_j)$ can be identified based on the thresholds (i.e., $d_j/(\aaaa_j^T \Sigma \aaaa_j + 1)^{1/2}, \forall j$) and the pairwise tetrachoric correlations (i.e., $(\aaaa_{j_1}^T \Sigma \aaaa_{j_2})/ \big((\aaaa_{j_1}^T \Sigma \aaaa_{j_1} + 1) (\aaaa_{j_2}^T \Sigma \aaaa_{j_2} + 1)\big)^{1/2}, \forall j_1\neq j_2$). 
The theoretical results on the sufficient conditions turned out to be very similar to those under the linear bifactor model and extensions. %The subsequent subsections present the identifiability results of each dichotomous response model (i.e., bifactor model, extended bifactor model, and two-tier item factor model) in details.

\subsection{Standard bifactor model}\label{pro:bifactor}

Adopting the same notations as the linear bifactor model, under the probit bifactor model, the probability of a response of $1$ on item $j\in \{1,\ldots, J\}$ in testlet $g_j$ is given by 
\begin{equation}\label{probit_bi}
P(Y_j = 1 | \eta_0, \eta_1, \ldots, \eta_G) = \Phi(d_j + a_0\eta_0 + \sum_{g=1}^G a_{jg}\eta_{g}) =\Phi(d_j + a_0\eta_0 + a_{jg_j}\eta_{g_j}),
\end{equation}
where, similar to the linear case, $a_{jg} = 0$ for all $g\neq g_j$, and $\eeta = (\eta_0, \eta_1, \ldots, \eta_G)^T\sim MVN(\mathbf 0, \Sigma)$ with $\Sigma = \boldsymbol I_{(1+G)\times (1+G)}.$ With $A$ and $\mathbf d$ denoting the loading matrix and the vector of intercepts, respectively, the parameter space of the probit standard bifactor model is given by $\mathfrak P = \{(A, \mathbf d) ~ |~ \textrm{first nonzero}$ $\textrm{element in every column of $A$ is positive} \}.$ Based on this, the definition of probit bifactor model identifiability is given below.
\begin{definition}
	We say a probit bifactor model is identifiable at $(A, \mathbf d)$ if for any other set of parameter $(A^{'}, \mathbf d^{'})$ that defines the same probability distribution, it must hold that 
	\begin{eqnarray}
	(A, \mathbf d) = (A^{'}, \mathbf d^{'}).
	\end{eqnarray} 
\end{definition}

Adopting the same definitions of sets $\H{1}, \H{2},$ and $\mathcal Q_g$ as in section \ref{lin:bifactor}, the theorem below provides the sufficient and necessary conditions for the identifiability of dichotomous bifactor models with probit link. 

\begin{theorem}\label{thm:probit:bifactor}
	Under standard bifactor model with probit link, the model parameter is identifiable if and only if it satisfies one of the follow conditions.
	\begin{itemize}
		\item[P1] $|\H{1}| \geq 3$; $|\mathcal Q_g| \geq 3$ for all $g = 1, \ldots, G$.
		\item[P2] $|\H{1}| = 2 $; $\H{2}$ is non-empty; $|\mathcal Q_g| \geq 3$ for $g = 1, \ldots, G$.
	\end{itemize}
\end{theorem}

The interpretations of $\mathit{P1}$ and $\mathit{P2}$ remain the same as for the linear bifactor model in Section \ref{lin:bifactor}. In the supplementary file, a few examples are provided to illustrate how the identifiability of the probit bifactor model can be checked.
%\begin{example}
%	Consider a probit bifactor model for $J = 9$ items, where the true parameters are given by 
%	\begin{eqnarray}
%	A = 
%	\begin{pmatrix}
%	a & b & 0 & 0 \\
%	a & b & 0 & 0 \\
%	a & b & 0 & 0 \\
%	a & 0 & c & 0 \\
%	a & 0 & c & 0 \\
%	a & 0 & c & 0 \\
%	a & 0 & 0 & d \\
%	a & 0 & 0 & d \\
%	a & 0 & 0 & d \\
%	\end{pmatrix},
%	\quad \mathbf d = 
%	\begin{pmatrix}
%	d_1 \\
%	d_2 \\
%	d_3 \\
%	d_4 \\
%	d_5 \\
%	d_6 \\
%	d_7 \\
%	d_8 \\
%	d_9 \\
%	\end{pmatrix},
%	\end{eqnarray}
%	with $a,b,c,d \neq 0$. The parameter is identifiable by checking that it satisfies Condition P1.
%\end{example}
%
%\begin{example}
%	Consider a probit bifactor model with $J = 8$ items and true parameters
%	\begin{eqnarray}
%	A = 
%	\begin{pmatrix}
%	a & b & 0 & 0 \\
%	a & b & 0 & 0 \\
%	a & b & 0 & 0 \\
%	a & 0 & c & 0 \\
%	a & 0 & c & 0 \\
%	a & 0 & c & 0 \\
%	a & 0 & 0 & d \\
%	a & 0 & 0 & d \\
%	\end{pmatrix},
%	\quad \mathbf d = 
%	\begin{pmatrix}
%	d_1 \\
%	d_2 \\
%	d_3 \\
%	d_4 \\
%	d_5 \\
%	d_6 \\
%	d_7 \\
%	d_8 \\
%	\end{pmatrix},
%	\end{eqnarray}
%	with $a,b,c,d \neq 0$. This setting is not identifiable by checking that it fails to satisfy  either Condition P1 or Condition P2.
%\end{example}

\subsection{Extended bifactor model}

With the same item response function as the standard bifactor model, the extended probit bifactor model relaxes the assumption of $\Sigma = \boldsymbol{I}_{(1+G)\times (1+G)}$ by allowing correlations among $\eta_1, \ldots, \eta_G$. The covariance matrix for $\eeta$, $\Sigma$, hence takes the form of $\Sigma = \begin{pmatrix}
1& \mathbf 0^T \\
\mathbf 0^T  & \boldsymbol\Sigma_G \\
\end{pmatrix}$, with $\Sigma_G$ positive definite with diagonal entries of $1$ and no further restriction on off-diagonal entries.
Under the probit extended bifactor model, the parameter space is given by
$\mathfrak P = \{(A, \mathbf d, \Sigma_G) ~ |~  \textrm{ first non-zero element in evey column of $A$ is positive}, \mathrm{diag}(\Sigma_G) = 1, \Sigma_G ~ \textrm{positive definite}. \}.$ And the definition of probit extended bifactor model identifiability is as follows.

\begin{definition}
	We say a probit extended bifactor model is identifiable at $(A, \mathbf d, \Sigma_G)$ if for any other set of parameters $(A^{'}, \mathbf d^{'}, \Sigma_G^{'})$ that define the same probability distribution,  
	\begin{eqnarray}
	(A, \mathbf d, \Sigma_G) = (A^{'}, \mathbf d^{'}, \Sigma_G^{'})
	\end{eqnarray} 
	must hold. 
\end{definition}

Again, it can be shown that the sufficient conditions and necessary conditions for linear extended bifactor model still hold when the responses become binary. 

\begin{theorem}\label{thm:extendprobit}
	Under the probit extended bifactor model, the model parameters are identifiable if one of the following set of requirements is satisfied:
	\begin{itemize}
		\item [E1S~] $|\mathcal Q_g| \geq 3$ for all $g = 1, \ldots, G$; $|\H{3}| \geq 3$.
		\item [E2S~] $|\mathcal Q_g| \geq 3$ for all $g = 1, \ldots, G$; $|\H{3}| = 2$; $|\H{2}| \geq 1$. 
	\end{itemize}
\end{theorem}

\begin{theorem}\label{thm:extendprobit:nec}
	Under the extended bifactor model with probit link, the model parameters are identifiable only  if both conditions below are satisfied.
	\begin{itemize}
		\item [E1N~] $|\mathcal Q_g| \geq 2$ for all $g = 1, \ldots, G$; and $|\mathcal Q_g| \geq 3$ for $g: \Sigma_G[g,-g] = \mathbf 0$. 
		\item [E2N~] $|\H{3}| \geq 2$. 
	\end{itemize}
\end{theorem}

Here, the definitions of the sets $\mathcal Q_g, \H{3},$ and $\H{2}$ and the interpretations of the conditions remain the same as those for the linear extended bifactor model in section \ref{lin:extended}.

\subsection{Two-tier model}
A two-tier probit model with $J$ items, $L$ main factors and $G$ testlets has the following item response function for a particular item $j$ in testlet $g_j$,
\begin{eqnarray}\label{probit:twotier}
P(Y_j = 1| \eeta) = \Phi(d_j + \sum_{l = 1}^L a_{jl} \eta_{l} + \sum_{g=L+1}^{L+G} a_{jg} \eta_{g}) =  \Phi(d_j + \sum_{l = 1}^L a_{jl} \eta_{l} +  a_{jg_j} \eta_{g_j}),
\end{eqnarray}
where $a_{jg} = 0, \forall g\neq g_j$. Same as the linear two-tier model, the latent traits $\eeta = (\eeta_1, \eeta_2)$, where $\boldsymbol \eta_1 = (\eta_1, \ldots, \eta_L)^T$ and $\boldsymbol \eta_2 = (\eta_{L+1}, \ldots, \eta_{L+G})^T$, are assumed to follow a multivariate normal distribution with mean $\mathbf 0$ and covariance matrix $\Sigma$, which takes the form of $\begin{pmatrix}
\Sigma_L& \mathbf 0^T \\
\mathbf 0^T  & I_{G \times G} \\
\end{pmatrix}$, with $\Sigma_L$ positive-definite with diagonal elements of 1 and off-diagonal elements between $-1$ and $1$. The parameter space for the probit two-tier model is hence given by 
$\mathfrak P = \{(A, \mathbf d, \Sigma_L) ~ |~ \textrm{first non-zero element in evey column of $A$ is positive}, \mathrm{diag}(\Sigma_L) = 1, \Sigma_L ~ \textrm{positive definite} \}$, and the definition of probit two-tier model identifiability is as follows.
\begin{definition}
	A probit two-tier model is identifiable at $(A, \mathbf d, \Sigma_L)$ if there is another set of parameters $(A^{'}, \mathbf d^{'}, \Sigma_L^{'})$ such that $\mathbf I(A) = \mathbf I(A^{'})$ and they define the same distribution,
	then it must hold that $(A, \mathbf d, \Sigma_L) = (A^{'}, \mathbf d^{'}, \Sigma_L^{'})$. 
\end{definition}

Below we provide a set of sufficient conditions for the identifiability of the probit two-tier model. 
\begin{theorem}\label{twotier:probit}
	Under the probit two-tier model, suppose true parameter satisfies Condition C1 and one of Conditions T1 - T3, then the parameter is identifiable:
	\begin{itemize}
		\item[T1S] $|\H{4}| \geq 3$, $A[\mathcal B_{\H{4}}, 1:L]$ contains an identity, where $\mathcal B_{\H{4}}$ is the set of items that make up the testlets in $\H{4}$.
		\item[T2S] $|\H{4}| \geq 2$, $|\H{5}| \geq 1$, $A[\mathcal B_{\H{4}}, 1:L]$ contains an identity.
		\item[T3S] A[:,1:L] contains an identity, and there exists a partition of testlets $ \{1, \ldots, G\} = \mathcal G_1 \cup \mathcal G_2$, such that (a) $\bar A_{\mathcal G_1}$ has full-column rank, and
		(b) $A[\mathcal B_{\mathcal G_2},1:L]$ has full-column rank. 
	\end{itemize}
\end{theorem}
Here, the definitions of the sets ($\mathcal H$s) and the interpretations of the conditions remain the same as for the linear two-tier model in section \ref{lin:twotier}.

\section{Remarks}

\vspace{-5mm}

\subsection{Orthogonality between primary and testlet dimensions}

\vspace{-5mm}

Discussions on the identification restriction for bifactor models can be found in \cite{rijmen2009efficient}, where it is pointed out that
three types of identification restrictions are required.
\begin{itemize}
	\item $G + 1$ restrictions for fixing the origins of general and testlet effects.
	\item $G + 1$ restrictions for fixing the scales of general and testlet effects.
	\item $G$ restrictions for dealing with the rotation issue. 
\end{itemize}

By translating the restrictions to mathematical expressions, the above three conditions are equivalent to 
\begin{eqnarray}
\eeta \sim N(\mathbf 0, \Sigma), \quad 
\Sigma = \begin{pmatrix}
1 & \mathbf 0^T\\
\mathbf 0 & \Sigma_G \\
\end{pmatrix},
\end{eqnarray}
with $\Sigma_G[g,g] = 1$ for all $g \in \{1, \ldots, G\}$.
However they do not provide a rigorous proof why we need the third type of restriction.
Below we provide Theorem \ref{bifactor:general:nonid} to answer this question.
Consider the parameter space 
\begin{eqnarray}\label{modelspace:3}
\mathfrak P = \{(A, \mathbf d, \Sigma_G, \rrho, \boldsymbol{\lambda}) ~ | ~  \mathrm{diag}(\Sigma) = 1, \Sigma ~\textrm{is positive definite}\}, 
% \Sigma_G[g,g^{'}] \in (-1,1) 
\end{eqnarray}
where
$\Sigma =  
\begin{pmatrix}
1 & \rrho^T \\
\rrho  & \Sigma_G \\
\end{pmatrix}.$
%\begin{eqnarray}\label{general:setting}
%\Sigma =  
%\begin{pmatrix}
%1 & \rrho^T \\
%\rrho  & \Sigma_G \\
%\end{pmatrix}.
%\end{eqnarray}
\begin{theorem}\label{bifactor:general:nonid}
	The bifactor model is not identifiable at any $(A, d, \Sigma_G, \rrho, \boldsymbol{\lambda})) \in \mathfrak P$ as defined in \eqref{modelspace:3}. 
\end{theorem}
The implication of Theorem \ref{bifactor:general:nonid} is that there is no identifiable model in $\mathfrak P$ as defined in \eqref{modelspace:3}, that is, when the orthogonality restriction between primary and testlet dimensions is further dropped. This explains why the identification results can only be extended to correlated testlet dimensions.

\subsection{Extensions}

The results here can be extended to more general settings. 
\begin{itemize}
	\item The normality assumptions in the linear bifactor model can be removed. That is, we do not require $\eta_g \sim N(0,1)$ and $\epsilon_j \sim N(0, \lambda_j)$ but instead assume $\textrm{Var}(\eta_g) = 1$ and $\textrm{Var}(\epsilon_j) = \lambda_j$. By checking the first and second moments, it is not hard to see that sufficient conditions in previous theorems still guarantee the identification.
	\item For the ordinal probit model, each $Y_j$ takes values in $\{1, \ldots, K_j\}$ $(K_j \geq 2)$ and follows the following probability distribution,
	\begin{eqnarray}
	P(Y_j > k\mid \eeta) = \Phi(d_j^{(k)} + \mathbf a_j^{T} \eeta)
	\end{eqnarray}
	for $k = 1, \ldots, K_j - 1$ with $d_j^{(1)} \geq d_j^{(2)} ... \geq d_j^{(K_j-1)}$. Under the same set of sufficient conditions, we can 
	easily obtain the identifiability results.
	
\end{itemize}

\subsection{Connections}

\vspace{-5mm}

Under the linear bifactor model setting, 
the sufficient condition in Theorem \ref{thm:generalother} given by \cite{anderson1956statistical} can be simplified in the sense that 
\textit{there are at least three items in each testlet, i.e. $|Q_g| \geq 3$ for all $g$}. 
(Suppose there exists a testlet with at most two items, then it is impossible to find two disjoint submatrices of $A$ with full column rank after deleting an item within that testlet.)
It can also be checked that this sufficient condition is satisfied by 
E1S and E2S in our Theorem \ref{thm:extendunknown}.  

For general linear factor models, two- and three- indicator rules are two sets of simple sufficient identifiability conditions; see \cite{bollen1989structural}, 
\begin{itemize}
	\item \textit{Two-indicator rules:}
	(1) Each latent factor is related to three items;
	(2)	Each row of $A$ has one and only one non-zero element;
	(3) Latent factors are uncorrelated.
	(4) $\epsilon$'s are uncorrelated.
	\item \textit{Three-indicator rules:}
	(1) Each latent factor is related to two items;
	(2)	Each row of $A$ has one and only one non-zero element;
	(3) No-zero elements in $\Sigma$;
	(4) $\epsilon$'s are uncorrelated.
\end{itemize} 

Although two- and three-indicator rules seem similar to the conditions in Theorems \ref{thm:unknown:bifactor} and \ref{thm:extendunknown}, they cannot be applied in bifactor/two-tier models.
By nature, it is impossible to assume each row of $A$ has one and only one non-zero element since each item has at least two latent dimensions (general factors and testlet-specific factor).
Fortunately, three items are enough for identifying the testlet effects. 
Thanks to the model structure, the general factor can also be identified when there is a sufficient number of testlets.

\section{Discussion}

This paper addresses the fundamental issue of identifiability of  bifactor model and its extensions, under both linear model with continuous indicators and probit model with dichotomous responses. The identifiability (or nonidentifiability) of a model can be determined through easily checkable conditions. 
In particular, conditions $\mathit {P1}$ and $\mathit {P2}$ establish the minimum requirements that can ensure the identifiability of the standard bifactor model. For the extended bifactor model with correlated subdimensions, a set of necessary conditions ($\mathit{E1N}, \mathit{E2N}$) and a set of sufficent conditions ($\mathit{E1S} - \mathit{E3S}$) for parameter identifiability were proposed. Sufficient conditions for two-tier model identifiability were further presented in $\mathit{C0}$ (or $\mathit {C1}$ for probit model) and $\mathit{T1S} - \mathit{T3S}$.
Theoretical results were able to explain underidentification phenomena observed in the existing literature. Simulation studies demonstrated the consequences on parameter estimation when the identifiability conditions were or were not met. From a practical viewpoint, these checkable identifiability conditions can guide test-developers through the design and evaluation of bifactor-type assessments.

It should be noted that, although both probit and logistic models can be applied for binary outcomes, the current identifiability results for probit models do not directly apply to item bifactor analysis with logistic parametrization, as seen in \cite{demars2006application}, \cite{cai2010two} and \cite{jeon2013modeling}. When a normal distribution is assumed for the latent traits, random-effect logistic item factor models involve the convolution of Gaussian and logistic random variables. This class of models hence do not imply the same first and second moments for item responses as the probit case. Future research may look into the identifiability conditions for bifactor-type models with logit link, perhaps adopting similar 
approaches as in \cite{san2013identification} for two-parameter logistic item response models. 
%Another natural extension of the current work is to probit models with polytomous responses. 
The current bifactor model identification findings may also be extended to higher-order factor models \cite{yung1999relationship}, under which latent factors are assumed to exhibit a hierarchical structure, with higher-order latent factors governing secondary, specific factors.

%%%%%%%%%%%%%%%%%%%%%%%%%%%%%%%%%%%%%%%%%%%%%%%%%%%%%%%%%%%%%%%%%%%%%%%%%%%%%%%%%%%%%%%%%%%%%%%%%%%%%%%%%%%%%%%%%%%%%%%%%%%%

%%%%%%%%%%%%%%%%%%%%%%%%%%%%%%%%%%%%%%%%%%%%%%%%%%%%%%%%%%%%%%%%%%%%%%%%%%%%%%%%%%%%%%%%%%%%%%%%%%%%%%%%%%%%%%%%%%%%%%%%%%%%
\section*{Supplementary Materials}

The supplementary material contains the simulation studies, illustrative examples and technical proofs of main theoretical results.

\par
%%%%%%%%%%%%%%%%%%%%%%%%%%%%%%%%%%%%%%%%%%%%%%%%%%%%%%%%%%%%%%%%%%%%%%%%%%%%%%%%%%%%%%%%%%%%%%%%%%%%%%%%%%%%%%%%%%%%%%%%%%%%
%\section*{Acknowledgements}
%
%Write the acknowledgements here.
%\par

%%%%%%%%%%%%%%%%%%%%%%%%%%%%%%%%%%%%%%%%%%%%%%%%%%%%%%%%%%%%%%%%%%%%%%%%%%%%%%%%%%%%%%%%%%

\bibhang=1.7pc
\bibsep=2pt
\fontsize{9}{14pt plus.8pt minus .6pt}\selectfont
\renewcommand\bibname{\large \bf References}
%\begin{thebibliography}{11}
\expandafter\ifx\csname
natexlab\endcsname\relax\def\natexlab#1{#1}\fi
\expandafter\ifx\csname url\endcsname\relax
  \def\url#1{\texttt{#1}}\fi
\expandafter\ifx\csname urlprefix\endcsname\relax\def\urlprefix{URL}\fi

%% use bibfile 
 \bibliographystyle{chicago}      % Chicago style, author-year citations
%  \bibliography{bibfile}   % name your BibTeX data base

%\newpage

%\bibliographystyle{plain}
\bibliography{references} 

%%  Another method

%%%%%%%%%%%%%%%%%%%%%%%%%%%%%%%%%%%%%%%%%%%%%%%%%%%%%%%%%%%%%%%%%%%%%%%%%%%%%%%%%%%%%%%%%%%%%%%%%%%%%%%%%%%%%%%%%%%%%%%%%%%%
\vskip .65cm
\noindent
Correspondence: Susu Zhang 
\vskip 2pt
\noindent
E-mail: szhan105@illinois.edu
\vskip 2pt

%\noindent
%second author affiliation
%\vskip 2pt
%\noindent
%E-mail: (second author email)

% \vskip .3cm
%\centerline{(Received ???? 20??; accepted ???? 20??)}\par

%%%%%%%%%%%%%%%%%%%%%%%%%%%%%%%%%%%%%%%%%%%%%%%5
\clearpage

\begin{center}
	{\large \textbf{Supplementary to ``Identifiabliity of the Bifactor Models"} }
\end{center}

	In this supplementary, we provide simulation results, additional examples and technical proofs for all theoretical results stated in the main paper.
	In particular, Appendix \ref{sec:simulation} provides simulation results for validating our theories.
	Appendix \ref{sec:example} gives multiple examples to help readers to better understand the structure of bifactor models.
	Appendix \ref{sec:C} collects the proofs for the identifiability of standard bifactor models.
	Appendix \ref{sec:D} is for the extended bifactor models.
	Appendix \ref{sec:E} gives the proofs of two-tier model's identifiabiltiy.
	Finally, the proofs of results in Section 4 can be found in Appendix \ref{sec:F}.
	\bigskip

	\section{Simulation}\label{sec:simulation}
	
	This section presents the results from several simulation studies designed to verify the theoretical identifiability results. Specifically, section \ref{sec:simbif} presents the numerical results on the probit bifactor model, and section \ref{sec:simebif} presents the results on the probit extended bifactor model. For both classes of models, several tests were considered: Some of the tests have true parameters that satisfy the identifiability conditions, while others fail to meet the conditions. For each test, $500$ sets of responses were randomly generated, with random samples of size $N = 1000, 2000,$ or $4000$.The stochastic expectation-maximization (StE) algorithm \citep{celeux1985sem,ip2002single} was employed to estimate both classes of models, which has been applied to IRT models \citep{diebolt1996stochastic,fox2003stochastic,zhang2020improved}. Specifically, a Gibbs sampler following \cite{albert1992bayesian} was adopted for the stochastic expectation (StE) step and a gradient descent algorithm for the maximization (M) step. To guarantee the convergence, the StE algorithm was iterated 10,000 times for each set and the first 5,000 iterations were discarded as burn-in. The estimated parameters ($\hat A, \hat{\mathbf{d}}$) were evaluated in terms of root mean squared error (RMSE) with respect to the true parameters. That is, for a particular entry of the $A$ matrix, $a_{jk},$ its RMSE was given by
	\begin{equation}
	RMSE(\hat a_{jk}) = (\frac{1}{500}\sum_{r=1}^{500} (\hat a_{jk} - a_{jk})^2)^{1/2}.
	\end{equation}
	
	For tests that meet the identifiability requirements, parameter estimates are expected to converge to the true values as $N$ increases, with the $RMSE$s approaching $0$. This will not be the case for tests that fail to meet the identifiability conditions, in which case the parameters cannot be consistently estimated.

	\subsection{Study 1: Probit bifactor model}\label{sec:simbif}
	
	Under the probit bifactor model, parameter recovery was evaluated under the following $4$ cases:
	\begin{itemize}
		\item[Case 1] Consider a test with three testlets. The true parameters are listed in Table \ref{Table:True}. By checking that $|\H{1}| = 3$ and $|\mathcal Q_g| > 3$ for all $g$, the model is identifiable according to Theorem \ref{thm:unknown:bifactor}.
		\item[Case 2] Remove testlet $3$ from Case 1. The model is no longer identifiable according to Theorem \ref{thm:unknown:bifactor}, because $|\H{1}| = 2$ and $|\H{2}| = 0$.
		\item[Case 3] Remove testlet $2$ from Case 1.  According to Theorem \ref{thm:unknown:bifactor}, the model is identifiable, because $|\H{1}| = 2$, $|\H{2}| = 1$ and $|\mathcal Q_g| > 3$ for all $g$.
		\item[Case 4] Based on the true parameters in Table \ref{Table:True}, construct a new test containing item 1 from testlet 1, items 11 and 12 from testlet 2, and all the 10 items in testlet 4. The model is nonidentifiable according to Theorem \ref{thm:unknown:bifactor} by checking that $|\mathcal Q_1|, |\mathcal Q_2 |< 3$.
	\end{itemize}

	\begin{table}
		\caption{True item parameters.}\label{Table:True}
		\centering
		{\scriptsize
			\begin{tabular}{rrrrrr}
				\hline
				&&\multicolumn{3}{c}{Testlet-specific factors} & \\
				item & Main factor & 1&2&3 &d   \\ 
				\hline
				1 & 1.00 & 2.00 &  &  & 1.51 \\ 
				2 & 1.00 & 2.00 &  &  & .39 \\ 
				3 & 1.00 & 2.00 &  &  & -.62 \\ 
				4 & 1.00 & 2.00 &  &  & -2.21 \\ 
				5 & 1.00 & 2.00 &  &  & 1.12 \\ 
				6 & 1.00 & 2.00 &  &  & -.04 \\ 
				7 & 1.00 & 2.00 &  &  & -.02 \\ 
				8 & 1.00 & 2.00 &  &  & .94 \\ 
				9 & 1.00 & 2.00 &  &  & .82 \\ 
				10 & 1.00 & 2.00 &  &  & .59 \\ 
				11 & 2.00 &  & 1.00 &  & .92 \\ 
				12 & 2.00 &  & 1.00 &  & .78 \\ 
				13 & 2.00 &  & 1.00 &  & .07 \\ 
				14 & 2.00 &  & 1.00 &  & -1.99 \\ 
				15 & 2.00 &  & 1.00 &  & .62 \\ 
				16 & 2.00 &  & 1.00 &  & -.06 \\ 
				17 & 2.00 &  & 1.00 &  & -.16 \\ 
				18 & 2.00 &  & 1.00 &  & -1.47 \\ 
				19 & 2.00 &  & 1.00 &  & -.48 \\ 
				20 & 2.00 &  & 1.00 &  & .42 \\ 
				21 & 1.00 &  &  & -.63 & 1.36 \\ 
				22 & 1.00 &  &  & .18 & -.10 \\ 
				23 & 1.00 &  &  & -.84 & .39 \\ 
				24 & 1.00 &  &  & 1.60 & -.05 \\ 
				25 & 1.00 &  &  & .33 & -1.38 \\ 
				26 & 1.00 &  &  & -.82 & -.41 \\ 
				27 & 1.00 &  &  & .49 & -.39 \\ 
				28 & 1.00 &  &  & .74 & -.06 \\ 
				29 & 1.00 &  &  & .58 & 1.10 \\ 
				30 & 1.00 &  &  & -.31 & .76 \\ 
				\hline
			\end{tabular}
		}
	\end{table}

	The average RMSEs for the $4$ cases, across all  non-zero $\hat a$s and $\hat d$s, are reported in Table \ref{Table:SIM_RMSE}.
	Compared to the two identifiable cases (Case 1 and Case 3), the RMSEs from the two unidentifiable cases (Case 2 and Case 4) were remarkably larger. This is consistent with the theoretical results on the sufficient and necessary condition for bifactor models. Moreover, under Case 4 where $|\mathcal Q_1|, |\mathcal Q_2 |< 3$ and $|\mathcal Q_3|>3$, the average RMSEs of $\hat a$s in testlet 1 were $1.65$, $.89$ and $.64$, respectively, for $N=1000, 2000$ and $4000$. For testlet 2, the average RMSEs of $\hat a$s were $1.33, 0.91$ and $0.58$, respectively. However, for testlet 3, the average RMSEs were $0.12, 0.07$ and $0.05$. This suggests that, for a particular testlet $g$, the parameters were better recovered when the requirement of $|\mathcal Q_g|\geq 3$ was met.
	
	\begin{table}
		\caption{RMSE of model parameters for bifactor models.}
		\label{Table:SIM_RMSE}
		\centering
		\begin{tabular}
			[c]{ccccccccccccccccc}\hline
			& \multicolumn{3}{c}{Case 1} & &  \multicolumn{3}{c}{Case 2}
			&  & \multicolumn{3}{c}{Case 3}&  & \multicolumn{3}{c}{Case 4}\\
			\cline{2-4}\cline{6-8}
			\cline{10-12}\cline{14-16}
			n	& 1000 & 2000 &4000 &&  1000 & 2000 &4000&  & 1000 & 2000&4000&& 1000 & 2000&4000\\\hline
			$\boldsymbol{a}$  & .16 & .10 & .07 && .77 & .64 & .55 && .17 & .11 & .07 && .49 & .30 & .20 \\
			$\mathbf d$& .11 & .07 & .05 && .18 & .10 & .06& & .11 & .07 & .05& & .30 & .18 & .12 \\ 
			\hline
		\end{tabular}
	\end{table}

	\subsection{Study 2: Probit extended bifactor model}\label{sec:simebif}
	Under the probit extended bifactor model, the following two cases were considered:
	
	\begin{itemize}
		\item[Case 5] 
		The three-testlet extended bifactor model with loadings specified in Table \ref{Table:True}.
		According to Theorem \ref{thm:extendunknown:nec}, we know the model is not identifiable by checking that $|\H{3}| = 1$. 
		\item[Case 6] 
		Replace testlet 2 in Table \ref{Table:True} by testlet 4 in Table \ref{Table:True_2}.
		According to Theorem \ref{thm:extendunknown}, the model is identifiable by observing that $|\H{3}| = 2$, $|\H{2}| = 2$ and $|\mathcal Q_g| \geq 3$ for $g = 1,3,4$.
	\end{itemize}
	
	For both cases, the true covariance matrix is given in Table \ref{Table:True_3}.
	
	Table \ref{Table:SIM_2_RMSE} reports the average RMSEs of the item parameters for the two cases with different sample sizes. Table \ref{Table:COV_2_RMSE_1000} provides the RMSE for each entry of $\Sigma_G,$ the covariance among the testlet factors.
	The numerical findings again corroborate the theoretical results, where both item and covariance parameters were recovered remarkably better under Case 6 compared to Case 5.
	
	\begin{table}[ht]
		\centering
		\caption{True parameters of testlet 4 for Case 6.}
		\label{Table:True_2}
		\begin{tabular}{rrrr}
			\hline
			&&Testlet-specific factors&\\
			item & Main factor & 4 &d   \\ 
			\hline
			31 & 2.00 & -.56 & -.16 \\ 
			32 & 2.00 & -.23 & -.25 \\ 
			33 & 2.00 & 1.56 & .70 \\ 
			34 & 2.00 & .07 & .56 \\ 
			35 & 2.00 & .13 & -.69 \\ 
			36 & 2.00 & 1.72 & -.71 \\ 
			37 & 2.00 & .46 & .36 \\ 
			38 & 2.00 & -1.27 & .77 \\ 
			39 & 2.00 & -.69 & -.11 \\ 
			40 & 2.00 & -.45 & .88 \\ 
			\hline
		\end{tabular}
	\end{table}
	
	\begin{table}[ht]
		\centering
		\caption{True covariance matrix for the extended bifactor model.}
		\label{Table:True_3}
		\begin{tabular}{rrrrrr}
			\hline
			&&\multicolumn{4}{c}{Testlet-specific factors}\\
			& Main factor & 1&2&3&4    \\ 
			\hline
			Main factor  & 1.00 &  &  & &     \\ 
			\multirow{4}{*}{Testlet-specific factors} & & 1.00 & & &  \\ 
			&  & .44 & 1.00 &  & \\ 
			&  & .32 & .52 & 1.00 &  \\ 
			&  & .26 & .21 & .29 &1.00   \\ 
			\hline
		\end{tabular}
	\end{table}

	\begin{table}
		\caption{RMSE of item parameters for extended bifactor model, under Case 5 and Case 6.}
		\label{Table:SIM_2_RMSE}
		\centering
		\begin{tabular}
			[c]{cccccccccc}\hline
			& \multicolumn{3}{c}{Case 5} & &  \multicolumn{3}{c}{Case 6}
			\\
			\cline{2-4}\cline{6-8}
			
			N	& 1000 & 2000 &4000 &&  1000 & 2000 &4000& \\\hline
			$\mathbf a$ & .46 & .37 & .32 && .20 & .11 & .07 \\ 
			d & .11 & .07 & .05& & .11 & .07 & .05 \\ 
			\hline
		\end{tabular}
	\end{table}
	
	\begin{table}
		\centering
		\caption{RMSE of the covariance matrix under Cases 5 and 6.}
		\label{Table:COV_2_RMSE_1000}
		\scalebox{0.6}{
			\begin{tabular}{ccccccc}
				\hline
				$N = 1000$ &\multicolumn{3}{c}{Case 5} &\multicolumn{3}{c}{Case 6}\\
				\hline
				\multirow{3}{*}{$\Sigma_G$}  & .00 &  &  & .00 &  &  \\ 
				& .11 & .00 &  & .04 & .00 & \\ 
				& .16 & .32 & .00 & .04 & .09 & .00  \\ 
				\hline
			\end{tabular}
			\bigskip
			\begin{tabular}{ccccccc}
				\hline
				$N = 2000$ &\multicolumn{3}{c}{Case 5} &\multicolumn{3}{c}{Case 6}\\
				\hline
				\multirow{3}{*}{$\Sigma_G$}  & .00 &  &  & .00 &  &  \\ 
				& .11 & .00 &  & .03 & .00 & \\ 
				& .17 & .30 & .00 & .03 & .04 & .00  \\ 
				\hline
			\end{tabular}
			\bigskip
			\begin{tabular}{ccccccc}
				\hline
				$N = 4000$ &\multicolumn{3}{c}{Case 5} &\multicolumn{3}{c}{Case 6}\\
				\hline
				\multirow{3}{*}{$\Sigma_G$}  & .00 &  &  & .00 &  &  \\ 
				& .12 & .00 &  & .02 & .00 & \\ 
				& .18 & .27 & .00 & .02 & .02 & .00  \\ 
				\hline
			\end{tabular}
		}
	\end{table}

	\newpage
	
	\section{Illustrative Examples}\label{sec:example}
	
	\begin{example}
		Consider a standard bifactor model for $J = 7$ items, where the true parameters are
		\begin{eqnarray}
		A = 
		\begin{pmatrix}
		a & b & 0 & 0 \\
		a & b & 0 & 0 \\
		c & d & 0 & 0 \\
		c & d & 0 & 0 \\
		e & 0 & f & 0 \\
		e & 0 & f & 0 \\
		e & 0 & f & 0 \\
		\end{pmatrix},
		\quad \mathbf d = 
		\begin{pmatrix}
		d_1 \\
		d_2 \\
		d_3 \\
		d_4 \\
		d_5 \\
		d_6 \\
		d_7 \\
		\end{pmatrix},
		\end{eqnarray}
		with $a,b,c,d,e,f \neq 0$ and $a d \neq b c$. One can check that (1) $|\H{1}| = 2$, because both testlet $1$ and testlet $2$ have nonzero main factor loadings; (2) $|\mathcal Q_g|\geq 3$ for all $g$, because $|\mathcal Q_1| = 4, |\mathcal Q_2| = 3$; and (3) $|\H{2}|\geq 1$, because testlet $1$ can be partitioned to $\mathcal B_{g,1} = \{1,3\}, \mathcal B_{g,2} = \{2,4\}$, each containing linearly independent columns. By Condition P2, we know that the model is identifiable.
	\end{example}

	\begin{example}
		Even though the main factors by themselves satisfy the identification conditions of factor models, in a two-tier model context, the main factor loadings can still be indistinguishable. Consider a three-testlet model with two main factors, where
		\begin{eqnarray}
		A = 
		\begin{pmatrix}
		1 & 0 & 2 & 0 & 0 \\
		1 & 0 & 3 & 0 & 0 \\
		1 & 0 & 4 & 0 & 0 \\
		0 & 1 & 0 & 1 & 0 \\
		0 & 1 & 0 & 1 & 0 \\
		0 & 1 & 0 & 1 & 0 \\
		0 & 1 & 0 & 0 & 3 \\
		0 & 1 & 0 & 0 & 3 \\
		0 & 1 & 0 & 0 & 3 \\
		\end{pmatrix},
		\quad 
		\Sigma_L = I_{3 \times 3}.
		\end{eqnarray}
		One can observe that the main factor loadings, $A_{:,1:2}$, satisfy the sufficient condition for factor model identifiability per the 3-indicator rule \citep[see][]{bollen1989structural}. In addition, there are three testlets with nonzero main factor loadings. However, we can easily construct another set of parameters with the same observed data distribution, say,
		\begin{eqnarray}
		A^{'} = 
		\begin{pmatrix}
		\sqrt{3/2} & 0 & \sqrt{1/2} & 0 & 0 \\
		\sqrt{4/2} & 0 & \sqrt{2/2} & 0 & 0 \\
		\sqrt{5/2} & 0 & \sqrt{3/2} & 0 & 0 \\
		0 & 1 & 0 & 1 & 0 \\
		0 & 1 & 0 & 1 & 0 \\
		0 & 1 & 0 & 1 & 0 \\
		0 & 1 & 0 & 0 & 3 \\
		0 & 1 & 0 & 0 & 3 \\
		0 & 1 & 0 & 0 & 3 \\
		\end{pmatrix},
		\quad 
		\Sigma_L^{'} = I_{2 \times 2}.
		\end{eqnarray}  
		
		This is because main factor 1 only depends on one testlet 1 and is consequently mixed up with the testlet-specific factor.
		
		We can easily construct another set of parameters.
		\begin{eqnarray}
		A^{''} = 
		\begin{pmatrix}
		1 & 0 & 2 & 0 & 0 \\
		1 & 0 & 3 & 0 & 0 \\
		1 & 0 & 4 & 0 & 0 \\
		0 & \sqrt{1/2} & 0 & \sqrt{3/2} & 0 \\
		0 & \sqrt{1/2} & 0 & \sqrt{3/2} & 0 \\
		0 & \sqrt{1/2} & 0 & \sqrt{3/2} & 0 \\
		0 & \sqrt{2} & 0 & 0 & 2\sqrt{2} \\
		0 & \sqrt{2} & 0 & 0 & 2\sqrt{2} \\
		0 & \sqrt{2} & 0 & 0 & 2\sqrt{2} \\
		\end{pmatrix},
		\quad 
		\Sigma_L^{''} = I_{2 \times 2}.
		\end{eqnarray} 
		This time, main factor 2 depends on two testlets and it is mixed up with second and third testlet-specific factors. 
	\end{example}

	\begin{example}
		It should further be noted that having testlets that load on multiple main factors would not suffice for two-tier model identifiability.
		Consider a two-tier model with three testlets and three main factors, where
		\begin{eqnarray}
		A = 
		\begin{pmatrix}
		0 & 1 & -1 & 1 & 0 & 0 \\
		0 & 1 & -1 & 2 & 0 & 0 \\
		0 & 1 & -1 & 3 & 0 & 0 \\
		2 & 0 & 1 & 0 & 3 & 0 \\
		2 & 0 & 1 & 0 & 2 & 0 \\
		2 & 0 & 1 & 0 & 1 & 0 \\
		1 & 1 & 0 & 0 & 0 & 2 \\
		1 & 1 & 0 & 0 & 0 & 3 \\
		1 & 1 & 0 & 0 & 0 & 1 \\
		\end{pmatrix},
		\quad 
		\Sigma_L = I_{3 \times 3}.
		\end{eqnarray}  
		
		Though each of the main factors is associated with multiple testlets, we can construct another set of parameters which implies the same distribution. 
		\begin{eqnarray}
		A^{'} = 
		\begin{pmatrix}
		0 & - \frac{3}{\sqrt{6}} & -\frac{1}{\sqrt{2}} & 1 & 0 & 0 \\
		0 & - \frac{3}{\sqrt{6}} & -\frac{1}{\sqrt{2}} & 2 & 0 & 0 \\
		0 & - \frac{3}{\sqrt{6}} & -\frac{1}{\sqrt{2}} & 3 & 0 & 0 \\
		\sqrt{3} & 0 & \sqrt{2} & 0 & 3 & 0 \\
		\sqrt{3} & 0 & \sqrt{2} & 0 & 2 & 0 \\
		\sqrt{3} & 0 & \sqrt{2} & 0 & 1 & 0 \\
		\frac{2}{\sqrt{3}} & -\frac{2}{\sqrt{6}} & 0 & 0 & 0 & 2 \\
		\frac{2}{\sqrt{3}} & -\frac{2}{\sqrt{6}} & 0 & 0 & 0 & 3 \\
		\frac{2}{\sqrt{3}} & -\frac{2}{\sqrt{6}} & 0 & 0 & 0 & 1 \\
		\end{pmatrix},
		\quad 
		\Sigma_L^{'} = I_{3 \times 3}.
		\end{eqnarray} 
	\end{example}

	\begin{example}
		Consider a probit bifactor model for $J = 9$ items, where the true parameters are given by 
		\begin{eqnarray}
		A = 
		\begin{pmatrix}
		a & b & 0 & 0 \\
		a & b & 0 & 0 \\
		a & b & 0 & 0 \\
		a & 0 & c & 0 \\
		a & 0 & c & 0 \\
		a & 0 & c & 0 \\
		a & 0 & 0 & d \\
		a & 0 & 0 & d \\
		a & 0 & 0 & d \\
		\end{pmatrix},
		\quad \mathbf d = 
		\begin{pmatrix}
		d_1 \\
		d_2 \\
		d_3 \\
		d_4 \\
		d_5 \\
		d_6 \\
		d_7 \\
		d_8 \\
		d_9 \\
		\end{pmatrix},
		\end{eqnarray}
		with $a,b,c,d \neq 0$. The parameter is identifiable by checking that it satisfies Condition P1.
	\end{example}
	
	\begin{example}
		Consider a probit bifactor model with $J = 8$ items and true parameters
		\begin{eqnarray}
		A = 
		\begin{pmatrix}
		a & b & 0 & 0 \\
		a & b & 0 & 0 \\
		a & b & 0 & 0 \\
		a & 0 & c & 0 \\
		a & 0 & c & 0 \\
		a & 0 & c & 0 \\
		a & 0 & 0 & d \\
		a & 0 & 0 & d \\
		\end{pmatrix},
		\quad \mathbf d = 
		\begin{pmatrix}
		d_1 \\
		d_2 \\
		d_3 \\
		d_4 \\
		d_5 \\
		d_6 \\
		d_7 \\
		d_8 \\
		\end{pmatrix},
		\end{eqnarray}
		with $a,b,c,d \neq 0$. This setting is not identifiable by checking that it fails to satisfy  either Condition P1 or Condition P2.
	\end{example}

	\begin{example}
		The theoretical results from the current paper provide explanations to the findings in existing studies on bifactor identification with rigor and generality.
		For example, \cite{green2018empirical} considered the empirical underidentification problem of bifactor model, which was encountered when fitting particular types of bifactor models to certain types of data sets.
		They demonstrated that the bifactor model can be underidentified in samples with homogenuous-within and homogenuous-between (HWHB) covariance structure, that is, 
		$\sigma_{j_1 j_2} = \sigma_{g_1 g_2}$ where $g_1$, $g_2$ are the testlets that items $j_1$, $j_2$ belong to.
		In particular, they considered the following loading matrices.
		\small{
			\begin{eqnarray}
			A = 
			\begin{pmatrix}
			0.8 & 0.3 & 0\\
			0.8 & 0.3 & 0\\
			0.8 & 0.3 & 0\\
			0.8 & 0.3 & 0\\
			0.8 & 0 & 0.3 \\
			0.8 & 0 & 0.3 \\
			0.8 & 0 & 0.3 \\
			0.8 & 0 & 0.3 \\
			\end{pmatrix},
			\quad
			A = 
			\begin{pmatrix}
			0.7& 0.4 & 0\\
			0.7& 0.4 & 0\\
			0.7& 0.4 & 0\\
			0.7& 0.4 & 0\\
			0.8 & 0 & 0.3 \\
			0.8 & 0 & 0.3 \\
			0.8 & 0 & 0.3 \\
			0.8 & 0 & 0.3 \\
			\end{pmatrix}.
			\end{eqnarray}
		}
		They showed that the above two bifactor models were not identifable by constructing different solutions which lead to the same model-implied covariance matrix.
		
		One can check that $|\H{1}| = 2$, $|\H{2}| = 0$ for both settings. It follows from the sufficient and necessary conditions in Theorem \ref{thm:unknown:bifactor} that the two models are not identifiable.
	\end{example}

	\newpage
	
	\section{Proofs for Standard Bifactor Models}\label{sec:C}
	
	\begin{proof}[Proof of Theorem \ref{thm:unknown:bifactor}]
		We first introduce a few more notations.
		\begin{itemize}
			\item Define $\mathcal Q_{0,g} = \{j ~ | ~ g_j = g, \mathbf a_0[j] \neq 0 \}$.
			\item Define $\H{6} = \{g ~ | ~ |\mathcal Q_{0,g}| \geq 2 \}$.
		\end{itemize}
		It is easy to see that $\H{2} \subset \H{6} \subset \H{1}$.
		
		We prove the results by considering the follow cases. We aim to show that the model is identifiable if and only if Case 1.a or Case 2.d holds.
		\begin{itemize}
			\item[] \textbf{Case 1}: $|\H{1}| \geq 3$. 
			\begin{itemize}
				\item[a]  $|\mathcal Q_g| \geq 3$ for all $g = 1, \ldots, G$.
				\item[b] $|\mathcal Q_g| \leq 2$ for some $g \in \{1, \ldots, G\}$.
			\end{itemize}  
			\item[] \textbf{Case 2}: $|\H{1}| = 2 $.
			\begin{itemize}
				\item[a] $|\mathcal Q_g| \leq 2$ for some $g \in \{1, \ldots, G\}$.
				\item[b] $|\mathcal Q_{0,g}| \leq 1$ for all $g \in \H{1}$, i.e., $\H{6}$ is empty.
				\item[c] 
				%			$|\mathcal Q_g| \geq 3$ for $g = 1, \ldots, G$;  $|\mathcal Q_{0,g}| \geq 2$ for some $g \in \H{1}$; 
				%			$|\H{3}| = 0$.
				$|\mathcal Q_g| \geq 3$ for $g = 1, \ldots, G$;
				$\H{6}$ is non-empty; $\H{2}$ is empty.	
				\item[{\color{black}d}] {\color{black} 
					$|\mathcal Q_g| \geq 3$ for $g = 1, \ldots, G$; $\H{2}$ is non-empty.}
			\end{itemize} 
			\item[] \textbf{Case 3}: $|\H{1}| \leq 1$.
		\end{itemize}
		
		First, we can see that if two sets of parameters lead to the same marginal distribution, it must hold that 
		\begin{eqnarray}\label{identity:unknown}
		A A^T + \Lambda = A^{'} (A^{'})^T + \Lambda^{'}
		\end{eqnarray}
		where $\Lambda = \mathrm{diag}((\lambda_1, \ldots \lambda_J))$. In other words, the off-diagonal elements are not collapsed with error variance.
		
		\textbf{Case 1} Suppose there is another set of parameters leading to the same model. Then we have ${\mathbf a}_0[\mathcal B_{g_1}] {\mathbf a}_0[\mathcal B_{g_2}]^T = {\mathbf a}_0^{'}[\mathcal B_{g_1}] {\mathbf a}_0^{'}[\mathcal B_{g_2}]^T$ for $g_1 \neq g_2 \in \H{1}$.
		Thus it implies that ${\mathbf a}_0[\mathcal B_{g}] = \pm {\mathbf a}_0^{'}[\mathcal B_{g}]$ for $g \in \H{1}$.
		By ${\mathbf a}_0[\mathcal B_{g_1}] {\mathbf a}_0[\mathcal B_{g}]^T = {\mathbf a}_0^{'}[\mathcal B_{g_1}] {\mathbf a}_0^{'}[\mathcal B_{g}]^T$ for $g_1 \in \H{1}$ and $g \notin \H{1}$, we further have ${\mathbf a}_0[\mathcal B_{g}] = \pm {\mathbf a}_0^{'}[\mathcal B_{g}]$ for $g \in \H{1}$.
		This implies that main factor loading is identifiable.
		
		Note the fact that ${\mathbf a}_g[j_1] {\mathbf a}_g[j_2] = {\mathbf a}_g^{'}[j_1] {\mathbf a}_g^{'}[j_2]$ for $j_1 \neq j_2 \in \mathcal B_g$ by comparing the correlations within testlet. By this, we consider the following.
		
		If $\mathcal Q_g \geq 3$, then we must have that ${\mathbf a}_g[j] = \pm {\mathbf a}_g^{'}[j]$ for $j \in \mathcal Q_g$. Further, it implies ${\mathbf a}_g[j] = \pm {\mathbf a}_g^{'}[j]$ for $j \notin \mathcal Q_g$. Thus Case 1.a is identifiable.
		
		If $\mathcal Q_g \leq 2$ for some $g$, we take out $j_1$ and $j_2$ from $\mathcal Q_g$ (if exist). Set ${\mathbf a}_g^{'}[j_1] = c \cdot {\mathbf a}_g[j_1]$ and ${\mathbf a}_g^{'}[j_2] = 1/c \cdot {\mathbf a}_g[j_2]$ for $c$ satisfying that 
		%$|(c^2 - 1| \cdot {\mathbf a}_g[j_1]^2 < \lambda_{j_1}$ \sz{remove $($?} and 
		\begin{eqnarray}\label{eqn:c} 
		|c^2 - 1| \cdot {\mathbf a}_g[j_1]^2 < \lambda_{j_1} ~ \textrm{and} ~ |1 - 1/c^2| \cdot {\mathbf a}_g[j_2]^2 < \lambda_{j_2}.
		\end{eqnarray}
		Such $c$ exists since that $c = 1$ is one of the solution. By continuity, we know that any $c$ sufficiently close to 1 satisfy \eqref{eqn:c} and keep the same sign of $\mathbf a_g$.
		This tells that $\mathbf a_g[j]$ is not uniquely determined. Hence Case 1.b is not identifiable.

		\textbf{Case 2}
		Suppose $\mathcal Q_g \leq 2$ for some $g$. By the same construction in Case 1.b, parameter $\mathbf a_g[j]$ cannot be identified for $j \in \mathcal Q_{g}$. Thus Case 2.a is not identifiable.
		
		Suppose $\mathcal Q_{0,g} \leq 1$ for $g \in \H{1}$. We can set ${\mathbf a}_0^{'}[j_1] = c \cdot {\mathbf a}_0[j_1]$ for $j_1 \in \mathcal Q_{0 g_1}$, $g_1 \in \H{1}$;
		set ${\mathbf a}_0^{'}[j_2] = 1 / c \cdot {\mathbf a}_0^{'}[j_2]$ for $j_2 \in \mathcal Q_{0 g_2}$, $g_2 \in \H{1}$
		and keep other $ a_g[j]$'s fixed. It is easy to check that ${\mathbf a}_g[j_1] {\mathbf a}_g[j_2] = {\mathbf a}_g^{'}[j_1] {\mathbf a}_g^{'}[j_2]$ holds for any $j_1 \neq j_2 \in \mathcal B_g$ and all $g$. We then choose $c$ close to 1 enough such that $|c^2 - 1| \cdot {\mathbf a}_0[j_1]^2 < \lambda_{j_1}$ and $|1 - 1/c^2| \cdot {\mathbf a}_0[j_2]^2 < \lambda_{j_2}$.  Then we can find $\lambda_{j_1}^{'}$ and $\lambda_{j_2}^{'}$ to keep \eqref{identity:unknown} hold, and the sign remains unchanged.
		Thus Case 2.b is not identifiable as the parameter $\mathbf a_0$ can not be determined.  
		
		%Suppose $|\mathcal Q_g| \geq 3$ for $g = 1, \ldots, G$;  $|\mathcal Q_{0,g}| \geq 2$ for some $g \in \H{1}$; $|\H{3}| = 0$.	
		%Then we know $\mathbf a_{g_1}[\mathcal B_{g_1}] = b_1 \cdot \mathbf a_0[\mathcal B_{g_1}]$ and $\mathbf a_{g_2}[\mathcal B_{g_2}] = b_2 \cdot \mathbf a_0[\mathcal B_{g_2}]$
		%for some $b_1, b_2 \neq 0$.
		%By comparing the off-diagonal elements within testlet, we then have 
		%\begin{eqnarray}\label{eqn:case:2c}
		%(1 - c^2 + b_1^2) {\mathbf a}_0[\mathcal B_{g_1} - \{j\}] {\mathbf a}_0[j] = 
		%{\mathbf a}_{g_1}^{'}[\mathcal B_{g_1} - \{j\}] {\mathbf a}_{g_1}^{'}[j]
		%\quad (j \in \mathcal B_{g_1})
		%\end{eqnarray}
		%and 
		%\begin{eqnarray}\label{eqn:case:2c-2}
		%(1 - 1/c^2 + b_2^2) {\mathbf a}_0[\mathcal B_{g_2} - \{j\}] {\mathbf a}_[j] = 
		%{\mathbf a}_{g_2}^{'}[\mathcal B_{g_2} - \{j\}] {\mathbf a}_{g_2}^{'}[j]
		%\quad (j \in \mathcal B_{g_2})
		%\end{eqnarray}
		%We can arbitrarily choose ${\mathbf a}_{g_1}^{'}[j] = \sqrt{(1 - c^2 + b_1^2)} \cdot {\mathbf a}_0[j]$ for $j \in \mathcal B_{g_1}$ and 
		%${\mathbf a}_{g_2}^{'}[j] = \sqrt{(1 - 1 / c^2 + b_2^2)} \cdot {\mathbf a}_0[j]$ for $j \in \mathcal B_{g_2}$ and $c \neq 1$. Thus $|\H{3}| \neq 0$.
		
		Suppose $|\mathcal Q_g| \geq 3$ for $g = 1, \ldots, G$; $|\mathcal Q_g| \geq 3$ for $g = 1, \ldots, G$; $\H{1} \cap \H{2} \cap \H{6}$ is non-empty.
		Let $g$ be the testlet satisfying that $g \in \H{1} \cap \H{2} \cap \H{6}$.		
		By comparing the off-diagonals of the covariance matrix, it must hold that 
		\begin{eqnarray}\label{eqn:case:3d:lin}
		(1 - c^2) {\mathbf a}_0[\mathcal B_{g} - \{j\}] {\mathbf a}_0[j] + {\mathbf a}_g[\mathcal B_{g} - \{j\}] {\mathbf a}_g[j]= 
		{\mathbf a}_{g}^{'}[\mathcal B_{g} - \{j\}] {\mathbf a}_{g}^{'}[j]
		\quad (\forall j \in \mathcal B_{g}).
		\end{eqnarray}
		By the property of $\H{2}$, we can find a partition of $\mathcal B_{g} = \mathcal B_{g,1} \cup \mathcal B_{g,2}$. Hence \eqref{eqn:case:3d:lin} can be written as
		\begin{eqnarray}\label{eqn:case:3d:mat:lin}
		(1 - c^2) {\mathbf a}_0[\mathcal B_{g,1}] {\mathbf a}_0[\mathcal B_{g,2}]^T + {\mathbf a}_g[\mathcal B_{g,1}] {\mathbf a}_g[\mathcal B_{g,2}]^T = 
		{\mathbf a}_{g}^{'}[\mathcal B_{g,1}] {\mathbf a}_{g}^{'}[\mathcal B_{g,2}]^T.
		\end{eqnarray}
		When $c^2 \neq 1$, the left hand side of \eqref{eqn:case:3d:mat:lin} has rank 2, while the right hand side of \eqref{eqn:case:3d:mat:lin} has at most rank 1.
		Thus $c^2 \equiv 1$, which reduces to Case 1.a. Hence Case 2.d is identifiable.

		Suppose $|\mathcal Q_g| \geq 3$ for $g = 1, \ldots, G$; $|\mathcal Q_g| \geq 3$ for $g = 1, \ldots, G$; $\H{1} \cap \H{2} \cap \H{6}$ is empty.
		By Case 2.b, we only need to consider the situation that $|\H{1} \cap \H{6}| \geq 1$ and $|\H{1} \cap \H{2} \cap \H{6} | = 0$.
		Take $g \in \H{1} \cap \H{6}$, we then know that 
		$\bar A_g$ can only take one of the following forms (after suitable row ordering),
		\begin{eqnarray}\label{eqn:Aform:lin}
		(1)
		\begin{pmatrix}
		a & b \\
		c & d \\
		e & f \\
		\end{pmatrix},
		\quad (2)
		\begin{pmatrix}
		a & b \\
		c & d \\
		c & d \\
		\mathbf c & \mathbf d \\
		\end{pmatrix},
		\end{eqnarray}
		where the matrix of form (1) is 3 by 2 and satisfies that $b,d,f \neq 0$ and at most one of $a,c,e$ is zero;
		the matrix of form (2) is $J_g$ ($J_g \geq 4$) by 2 and contains at most two different rows (i.e. rows are not equal up to scaling).
		Under both cases, we only need to check \eqref{eqn:case:3d:lin} for items corresponding to the first three rows. For notational simplicity, we denote three items as 1,2 and 3.  
		
		We can construct another set of parameters, where
		\begin{eqnarray}\label{construct:2c:lin}
		{\mathbf a}_g^{'}[1] &=& \frac{(1 - c^2) {\mathbf a}_0^{}[1] {\mathbf a}_0^{}[3] + {\mathbf a}_g^{}[1] {\mathbf a}_g^{}[3]}{{\mathbf a}_g^{}[3] x}; \nonumber \\
		{\mathbf a}_g^{'}[2] &=& \frac{(1 - c^2) {\mathbf a}_0^{}[2] {\mathbf a}_0^{}[3] + {\mathbf a}_g^{}[2] {\mathbf a}_g^{}[3]}{{\mathbf a}_g^{}[3] x}; \nonumber \\
		{\mathbf a}_g^{'}[3] &=& {\mathbf a}_g^{}[3] x; \nonumber \\
		x &=&  \frac{((1 - c^2) {\mathbf a}_0^{}[1] {\mathbf a}_0^{}[3] + {\mathbf a}_g^{}[1] {\mathbf a}_g^{}[3]) ((1 - c^2) {\mathbf a}_0^{}[2] {\mathbf a}_0^{}[3] + {\mathbf a}_g^{}[2] {\mathbf a}_g^{}[3])}{((1 - c^2) {\mathbf a}_0^{}[1] {\mathbf a}_0^{}[2] + {\mathbf a}_g^{}[1] {\mathbf a}_g^{}[2]) {\mathbf a}_g[3]^2}; \nonumber \\
		\lambda^{'}_{j} &=& \lambda_j + (1 - c^2) a_0^2[j] + (a_g[j])^2 - (a_g^{'}[j])^2, ~ j = 1,2,3.
		\end{eqnarray} 
		Notice that ${\mathbf a}_g[3] \neq 0$ according to the assumption that $|\mathcal Q_g| \geq 3$. Hence the above solution will be different from the true parameters when $c^2 \neq 1$. Note that $x > 0$ when $c$ is sufficiently close to $1$. Then $\mathbf a_g^{'}$ has the same sign as $\mathbf a_g$. Thus Case 2.c is not identifiable.
		
		\textbf{Case 3} Apparently, it is not identifiable. 
		This is because we can construct another set of parameters, 
		$\bar A_g^{'} = \bar A_g, g = 2, \ldots, G$ and $\bar A_1^{'} = \bar A_1 R$ with $R$ being a 2 by 2 rotation matrix. It is easy to see that the two sets of parameters lead to the same distribution, since $\Sigma_{g_1 g_2} = \Sigma_{g_1 g_2}^{'}$ for all $g_1$ and $g_2$. 
		In addition, we can easily choose the rotation matrix $R$ such that it keeps sign of first non-zero element in each column of $\bar A_1$.
		
		Identifiability of $\mathbf d$ is obvious by using the expectation of $Y_j$. Thus we conclude the proof.
	\end{proof}

	\bigskip
	
	\begin{proof}[Proof of Proposition \ref{prop:tetra}]
		The sufficient part is straightforward by noticing that $P(Y_{j_1} = 1, \ldots, Y_{j_k} = 1)$ only depends on $d_j / (\aaaa_j^T \Sigma \aaaa_j + 1)^{1/2}$'s and $ (\aaaa_{j_1}^T \Sigma \aaaa_{j_2}) /\big((\aaaa_{j_1}^T \Sigma \aaaa_{j_1} + 1)(\aaaa_{j_2}^T \Sigma \aaaa_{j_2} + 1)\big)^{1/2}$'s for all possible combinations of $j_1, \ldots  j_k$.
		
		The necessary part is also not hard. Notice that CDF function $\Phi$ is a strictly monotone increasing function. By \eqref{1dim}, we must have  
		$d_j/(\aaaa_j^T \Sigma \aaaa_j + 1)^{1/2} = d_j^{'}/((\aaaa_j^{'})^T \Sigma^{'} \aaaa_j^{'} + 1)^{1/2}$ for all $j$.
		In addition, $\Phi_2(a,b,\rho)$ is a strictly monotone increasing function of $\rho$ for any fixed $a, b$. Thus, from \eqref{2dim}, we get 
		$(\aaaa_{j_1}^T \Sigma \aaaa_{j_2})/\big((\aaaa_{j_1}^T \Sigma \aaaa_{j_1} + 1) (\aaaa_{j_2}^T \Sigma \aaaa_{j_2} + 1)\big)^{1/2} = 
		(\aaaa_{j_1}^T \Sigma^{'} \aaaa_{j_2})/\big(((\aaaa_{j_1}^{'})^T \Sigma^{'} \aaaa_{j_1}^{'} + 1) ((\aaaa_{j_2}^{'})^T \Sigma^{'} \aaaa_{j_2}^{'} + 1)\big)^{1/2}$ for any $j_1 \neq j_2$. Hence we prove the proposition.
	\end{proof}
	
	\bigskip
	
	\begin{proof}[Proof of Theorem \ref{thm:probit:bifactor}]
		We keep using the notations of $\mathcal Q_{0,g}$ and $\H{6}$. 
		We still aim to show that the model is identifiable if and only if Case 1.a or Case 2.d holds.
		\begin{itemize}
			\item[] \textbf{Case 1}: $|\H{1}| \geq 3$. 
			\begin{itemize}
				\item[{\color{black}a}] {\color{black} $|\mathcal Q_g| \geq 3$ for all $g = 1, \ldots, G$.}
				\item[b] $|\mathcal Q_g| \leq 2$ for some $g \in \{1, \ldots, G\}$.
			\end{itemize}  
			\item[] \textbf{Case 2}: $|\H{1}| = 2 $.
			\begin{itemize}
				\item[a] $|\mathcal Q_g| \leq 2$ for some $g \in \{1, \ldots, G\}$.
				\item[b] $|\mathcal Q_{0,g}| \leq 1$ for all $g \in \H{1}$, i.e., $\H{6}$ is empty.
				\item[c] 
				%			$|\mathcal Q_g| \geq 3$ for $g = 1, \ldots, G$;  $|\mathcal Q_{0,g}| \geq 2$ for some $g \in \H{1}$; 
				%			$|\H{3}| = 0$.
				$|\mathcal Q_g| \geq 3$ for $g = 1, \ldots, G$; $\H{2}$ is empty.	
				\item[{\color{black}d}] {\color{black} 
					$|\mathcal Q_g| \geq 3$ for $g = 1, \ldots, G$; $\H{2}$ is non-empty.}
			\end{itemize} 
			\item[] \textbf{Case 3}: $|\H{1}| \leq 1$.
		\end{itemize}
		
		The first step is to show that the mapping $(x,y) \rightarrow (\frac{x}{\sqrt{x^2 + y^2 + 1}}, \frac{y}{\sqrt{x^2 + y^2 + 1}})$ is one-to-one.
		It is easy to see that the mapping is onto. We only need to show it is injective. 
		If there exists another $(x^{'}, y^{'})$ such that $\frac{x}{\sqrt{x^2 + y^2 + 1}} = \frac{x^{'}}{\sqrt{(x^{'})^2 + (y^{'})^2 + 1}}$ and 
		$\frac{y}{\sqrt{x^2 + y^2 + 1}} = \frac{y^{'}}{\sqrt{(x^{'})^2 + (y^{'})^2 + 1}}$, we can find $(x^{'})^2 + (y^{'})^2 = x^2 + y^2$ which further implies $x = x^{'}$ and $y = y^{'}$.
		Therefore, in the following, we only need to work with $\tilde {\mathbf a}_{j}$, where $\tilde {\mathbf a}_{j} = \frac{\mathbf a_j}{\sqrt{\mathbf a_j^T \mathbf a_j + 1}}$.

		\textbf{Case 1} Suppose there is another set of parameters leading to the same distribution. Then we have $\tilde{\mathbf a}_0[\mathcal B_{g_1}] \tilde{\mathbf a}_0[\mathcal B_{g_2}]^T = \tilde{\mathbf a}_0^{'}[\mathcal B_{g_1}] \tilde{\mathbf a}_0^{'}[\mathcal B_{g_2}]^T$ for $g_1 \neq g_2 \in \H{1}$.
		Thus it implies that $\tilde{\mathbf a}_0[\mathcal B_{g}] = \pm \tilde{\mathbf a}_0^{'}[\mathcal B_{g}]$ for $g \in \H{1}$.
		By $\tilde{\mathbf a}_0[\mathcal B_{g_1}] \tilde{\mathbf a}_0[\mathcal B_{g}]^T = \tilde{\mathbf a}_0^{'}[\mathcal B_{g_1}] \tilde{\mathbf a}_0^{'}[\mathcal B_{g}]^T$ for $g_1 \in \H{1}$ and $g \notin \H{1}$, we further have $\tilde{\mathbf a}_0[\mathcal B_{g}] = \pm \tilde{\mathbf a}_0^{'}[\mathcal B_{g}]$ for $g \in \H{1}$.
		This implies that main factor loading is identifiable.
		
		Notice that $\tilde{\mathbf a}_g[j_1] \tilde{\mathbf a}_g[j_2] = \tilde{\mathbf a}_g^{'}[j_1] \tilde{\mathbf a}_g^{'}[j_2]$ for $j_1 \neq j_2 \in \mathcal B_g$ by comparing the correlations within testlet.
		
		If $\mathcal Q_g \geq 3$, then we must have that $\tilde{\mathbf a}_g[j] = \pm \tilde{\mathbf a}_g^{'}[j]$ for $j \in \mathcal Q_g$. Further, it implies $\tilde{\mathbf a}_g[j] = \pm \tilde{\mathbf a}_g^{'}[j]$ for $j \notin \mathcal Q_g$. Thus Case 1.a is identifiable.
		
		If $\mathcal Q_g \leq 2$ for some $g$, we take out $j_1$ and $j_2$ from $\mathcal Q_g$ (if exist). Set $\tilde{\mathbf a}_g^{'}[j_1] = c \cdot \tilde{\mathbf a}_g[j_1]$ and $\tilde{\mathbf a}_g^{'}[j_2] = 1/c \cdot \tilde{\mathbf a}_g[j_2]$ for $c$ satisfying that $ \tilde{\mathbf a}_0[j_1]^2 + c^2 \cdot \tilde{\mathbf a}_g[j_1]^2 < 1$ and 
		\begin{eqnarray}\label{eqn:c:probit} 
		\tilde{\mathbf a}_0[j_2]^2 + 1 / c^2 \cdot \tilde{\mathbf a}_g[j_2]^2 < 1.
		\end{eqnarray}
		Such $c$ exists since that $c = 1$ is one of the solution. By continuity, we know that any $c$ sufficiently close to 1 satisfy \eqref{eqn:c:probit} and keep the same sign of $\tilde{\mathbf a_g}$.
		This tells that $\mathbf a_g[j]$ is not uniquely determined. Hence Case 1.b is not identifiable.

		\textbf{Case 2}
		Suppose $\mathcal Q_g \leq$ for some $g$. By the same construction in Case 1.b, parameter $\mathbf a_g[j]$ cannot be identified for $j \in \mathcal Q_{g}$. Thus Case 2.a is not identifiable.
		
		Suppose $\mathcal Q_{0,g} \leq 1$ for $g \in \H{1}$. We can set $\tilde{\mathbf a}_0^{'}[j_1] = c \cdot \tilde{\mathbf a}_0[j_1]$ for $j_1 \in \mathcal Q_{0 g_1}$, $g_1 \in \H{1}$;
		set $\tilde{\mathbf a}_0^{'}[j_2] = 1 / c \cdot \tilde{\mathbf a}_0^{'}[j_2]$ for $j_2 \in \mathcal Q_{0 g_2}$, $g_2 \in \H{1}$
		and keep other $\tilde a_g[j]$'s fixed. It is easy to check that $\tilde{\mathbf a}_g[j_1] \tilde{\mathbf a}_g[j_2] = \tilde{\mathbf a}_g^{'}[j_1] \tilde{\mathbf a}_g^{'}[j_2]^T$ holds for any $j_1 \neq j_2 \in \mathcal B_g$ and all $g$. Here $c$ is chosen to be positive to keep the sign.
		Thus Case 2.b is not identifiable as the parameter $\mathbf a_0[j_1]$ can not be determined.  
		
		%Suppose $|\mathcal Q_g| \geq 3$ for $g = 1, \ldots, G$;  $|\mathcal Q_{0,g}| \geq 2$ for some $g \in \H{1}$; $|\H{3}| = 0$.	
		%Then we know $\mathbf a_{g_1}[\mathcal B_{g_1}] = b_1 \cdot \mathbf a_0[\mathcal B_{g_1}]$ and $\mathbf a_{g_2}[\mathcal B_{g_2}] = b_2 \cdot \mathbf a_0[\mathcal B_{g_2}]$
		%for some $b_1, b_2 \neq 0$.
		%By comparing the correlation within testlet, we then have 
		%\begin{eqnarray}\label{eqn:case:2c}
		%(1 - c^2 + b_1^2) \tilde{\mathbf a}_0[\mathcal B_{g_1} - \{j\}] \tilde{\mathbf a}_0[j] = 
		%\tilde{\mathbf a}_{g_1}^{'}[\mathcal B_{g_1} - \{j\}] \tilde{\mathbf a}_{g_1}^{'}[j]
		%\quad (j \in \mathcal B_{g_1})
		%\end{eqnarray}
		%and 
		%\begin{eqnarray}\label{eqn:case:2c-2}
		%(1 - 1/c^2 + b_2^2) \tilde{\mathbf a}_0[\mathcal B_{g_2} - \{j\}] \tilde{\mathbf a}_[j] = 
		%\tilde{\mathbf a}_{g_2}^{'}[\mathcal B_{g_2} - \{j\}] \tilde{\mathbf a}_{g_2}^{'}[j]
		%\quad (j \in \mathcal B_{g_2})
		%\end{eqnarray}
		%We can arbitrarily choose $\tilde{\mathbf a}_{g_1}^{'}[j] = \sqrt{(1 - c^2 + b_1^2)} \cdot \tilde{\mathbf a}_0[j]$ for $j \in \mathcal B_{g_1}$ and 
		%$\tilde{\mathbf a}_{g_2}^{'}[j] = \sqrt{(1 - 1 / c^2 + b_2^2)} \cdot \tilde{\mathbf a}_0[j]$ for $j \in \mathcal B_{g_2}$ and $c \neq 1$.
		%Thus Case 2.c is not identifiable.
		
		Suppose $|\mathcal Q_g| \geq 3$ for $g = 1, \ldots, G$; $|\mathcal Q_g| \geq 3$ for $g = 1, \ldots, G$; $\H{1} \cap \H{2} \cap \H{6}$ is non-empty.
		Let $g$ be the testlet satisfying that $g \in \H{1} \cap \H{2} \cap \H{6}$.		
		By comparing the correlation within matrix, it must hold that 
		\begin{eqnarray}\label{eqn:case:3d}
		(1 - c^2) \tilde{\mathbf a}_0[\mathcal B_{g} - \{j\}] \tilde{\mathbf a}_0[j] + \tilde{\mathbf a}_g[\mathcal B_{g} - \{j\}] \tilde{\mathbf a}_g[j]= 
		\tilde{\mathbf a}_{g}^{'}[\mathcal B_{g} - \{j\}] \tilde{\mathbf a}_{g}^{'}[j]
		\quad (\forall j \in \mathcal B_{g}).
		\end{eqnarray}
		By the property of $\H{2}$, we can find a partition of $\mathcal B_{g} = \mathcal B_{g,1} \cup \mathcal B_{g,2}$. Hence \eqref{eqn:case:3d} can be written as
		\begin{eqnarray}\label{eqn:case:3d:mat}
		(1 - c^2) \tilde{\mathbf a}_0[\mathcal B_{g,1}] \tilde{\mathbf a}_0[\mathcal B_{g,2}]^T + \tilde{\mathbf a}_g[\mathcal B_{g,1}] \tilde{\mathbf a}_g[\mathcal B_{g,2}]^T = 
		\tilde{\mathbf a}_{g}^{'}[\mathcal B_{g,1}] \tilde{\mathbf a}_{g}^{'}[\mathcal B_{g,2}]^T.
		\end{eqnarray}
		When $c^2 \neq 1$, the left hand side of \eqref{eqn:case:3d:mat} has rank 2, while the right hand side of \eqref{eqn:case:3d:mat} has at most rank 1.
		Thus $c^2 \equiv 1$, which reduces to Case 1.a. Hence Case 2.d is identifiable.

		Suppose $|\mathcal Q_g| \geq 3$ for $g = 1, \ldots, G$; $|\mathcal Q_g| \geq 3$ for $g = 1, \ldots, G$; $\H{1} \cap \H{2} \cap \H{6}$ is empty.
		By Case 2.b, we only need to consider the situation that $|\H{1} \cap \H{6}| \geq 1$ and $|\H{1} \cap \H{2} \cap \H{6} | = 0$.
		Take $g \in \H{1} \cap \H{6}$, we then know that 
		$\bar A_g$ can only take one of the following forms (after suitable row ordering),
		\begin{eqnarray}\label{eqn:Aform}
		(1)
		\begin{pmatrix}
		a & b \\
		c & d \\
		e & f \\
		\end{pmatrix},
		\quad (2)
		\begin{pmatrix}
		a & b \\
		c & d \\
		c & d \\
		\mathbf c & \mathbf d \\
		\end{pmatrix},
		\end{eqnarray}
		where the matrix of form (1) is 3 by 2 and satisfies that $b,d,f \neq 0$ and at most one of $a,c,e$ is zero;
		the matrix of form (2) is $J_g$ ($J_g \geq 4$) by 2 and contains at most two different rows (i.e. rows are not equal up to scaling).
		Under both cases, we only need to check \eqref{eqn:case:3d} for items corresponding to the first three rows. For notational simplicity, we denote three items as 1,2 and 3.  
		
		We can construct another set of parameters,
		\begin{eqnarray}\label{construct:2c}
		\tilde{\mathbf a}_g^{'}[1] &=& \frac{(1 - c^2) \tilde{\mathbf a}_0^{}[1] \tilde{\mathbf a}_0^{}[3] + \tilde{\mathbf a}_g^{}[1] \tilde{\mathbf a}_g^{}[3]}{\tilde{\mathbf a}_g^{}[3] x}; \nonumber \\
		\tilde{\mathbf a}_g^{'}[2] &=& \frac{(1 - c^2) \tilde{\mathbf a}_0^{}[2] \tilde{\mathbf a}_0^{}[3] + \tilde{\mathbf a}_g^{}[2] \tilde{\mathbf a}_g^{}[3]}{\tilde{\mathbf a}_g^{}[3] x}; \nonumber \\
		\tilde{\mathbf a}_g^{'}[3] &=& \tilde{\mathbf a}_g^{}[3] x; \nonumber \\
		x &=&  \frac{((1 - c^2) \tilde{\mathbf a}_0^{}[1] \tilde{\mathbf a}_0^{}[3] + \tilde{\mathbf a}_g^{}[1] \tilde{\mathbf a}_g^{}[3]) ((1 - c^2) \tilde{\mathbf a}_0^{}[2] \tilde{\mathbf a}_0^{}[3] + \tilde{\mathbf a}_g^{}[2] \tilde{\mathbf a}_g^{}[3])}{((1 - c^2) \tilde{\mathbf a}_0^{}[1] \tilde{\mathbf a}_0^{}[2] + \tilde{\mathbf a}_g^{}[1] \tilde{\mathbf a}_g^{}[2]) \tilde{\mathbf a}_g[3]^2}.
		\end{eqnarray} 
		Notice that $\tilde{\mathbf a}_g[3] \neq 0$ according to assumption that $|\mathcal Q_g| \geq 3$. Hence the above solution will be different from true parameter when $c^2 \neq 1$. When $c$ is close enough to $1$, we know that $x$ is positive and $\tilde{\mathbf a}_g^{'}$ has the same sign as that of $\tilde{\mathbf a}_g$. Thus Case 2.c is not identifiable.
		
		\textbf{Case 3} Obviously, it is not identifiable by the same reason as stated in Case 3 in the proof of Theorem \ref{thm:unknown:bifactor}.
		
		Once loading matrix is identifiable, we can immediately identify $\mathbf d$ by using \eqref{1dim}. Thus we conclude the proof.
		%----------------------------------------	
		%	let $j$ be the row with two non-zero elements. We then find that $\tilde{\mathbf a}_0[\mathcal B_{g} - \{j\}]$ and  $\tilde{\mathbf a}_g[\mathcal B_{g} - \{j\}]$ are linearly independent. Thus the left hand side of \eqref{eqn:case:3d} is of rank  	
		% ---------------------------------
		%	Let $a = \frac{x}{\sqrt{x^2 + y^2 + 1}}, b = \frac{y}{\sqrt{x^2 + y^2 + 1}}$. We can solve that $x = \frac{a}{\sqrt{1 - (a^2 + b^2)}}$ and 
		%	$y = \frac{b}{\sqrt{1 - (a^2 + b^2)}}$.
	\end{proof}
	
	\newpage

	\section{Proofs for Extended Bifactor Models}\label{sec:D}
	
	For the linear and probit extended bifactor model identifiability, we only provide the complete proof for the linear case (i.e., Theorem \ref{thm:extendunknown} and Theorem \ref{thm:extendunknown:nec}). It should be apparent from Appendix A that the proofs for linear and probit models are very similar. 
	
	%\sz{personally I think it is better to present the sufficiency proof for the linear case. Theorems \ref{extendedH2} and \ref{extended:suff:general} seem to be corresponding to the linear case. Proof for the necessary conditions is also for the linear case. }

	%\sz{I moved the theorem for condition $C1$ here, because the extended bifactor probit identifiability seems to use it. If this section no longer uses it, feel free to move it to Appendix C.}

	The following are two support theorems for the proof of Theorems \ref{thm:extendunknown} and \ref{thm:extendunknown:nec}.
	\begin{theorem}\label{extended:suff:general}
		Under the the linear extended bifactor model with known error variance, if parameters satisfy $|\H{3}| \geq 2$ and $|\mathcal N| = 0$, then it is identifiable. 
	\end{theorem}
	
	\begin{theorem}\label{extendedH2}
		Under the linear extended bifactor model, the model parameter is not identifiable if $|\H{3}| \leq 1$.
	\end{theorem}
	
	Before proof of Theorems \ref{extended:suff:general} and \ref{extendedH2}, we first state the following Lemma \ref{fund:lem} which plays an important role in proving the identification of the extended bifactor model.
	\begin{lemma}\label{fund:lem}
		Assume $A$ and $B$ are both two-column matrices. Let $\Sigma$ and $\Sigma^{'}$ be two by two matrices. Consider the following situations:
		\begin{enumerate}
			\item Suppose both $A$ and $B$ are full-column rank. Thus $A \Sigma B^T = A \Sigma^{'} B^T$ implies that $\Sigma = \Sigma^{'}$.
			\item Suppose $A$ is full-column rank and $B$ has column-rank 1, i.e., $B = \mathbf b \begin{pmatrix}b_1 \\ b_2
			\end{pmatrix}^T$. Thus $A \Sigma B^T = A \Sigma^{'} B^T$ implies that 
			$\Sigma 
			\begin{pmatrix} 
			b_1 \\
			b_2 \\
			\end{pmatrix} = 
			\Sigma^{'}
			\begin{pmatrix} 
			b_1 \\
			b_2 \\
			\end{pmatrix}.$
			\item 
			Suppose both $A$ and $B$ are column-rank 1, i.e., $A = \mathbf a \begin{pmatrix}a_1 \\ a_2
			\end{pmatrix}^T$ and $B = \mathbf b \begin{pmatrix}b_1 \\ b_2
			\end{pmatrix}^T$. Thus $A \Sigma B^T = A \Sigma^{'} B^T$ implies that 
			$
			\begin{pmatrix}a_1 \\ a_2
			\end{pmatrix}^T
			\Sigma 
			\begin{pmatrix} 
			b_1 \\
			b_2 \\
			\end{pmatrix} = 
			\begin{pmatrix}a_1 \\ a_2
			\end{pmatrix}^T
			\Sigma^{'}
			\begin{pmatrix} 
			b_1 \\
			b_2 \\
			\end{pmatrix}.$
		\end{enumerate}
	\end{lemma}
	
	\bigskip
	
	\begin{proof}[Proof of Theorem \ref{extended:suff:general}]
		We prove this by contradiction. Suppose there exists another set of $\{A_g^{'}\}$ leading to the same distribution. 
		We pick any item pair $g_1$ and $g_2$ from set $\mathcal H^{\ast}$. We know that  $\Sigma_{gg} = \Sigma_{gg}^{'}$ for $g = g_1, g_2$, which implies that 
		\begin{eqnarray}
		\bar A_{g_1}^{'} = \bar A_{g_1} 
		\begin{pmatrix} 
		\cos \theta_1 & \sin \theta_1  \\
		- \sin \theta_1 & \cos \theta_1  \\
		\end{pmatrix}
		, 
		\qquad
		\bar A_{g_2}^{'} = \bar A_{g_2} 
		\begin{pmatrix} 
		\cos \theta_2 & \sin \theta_2  \\
		- \sin \theta_2 & \cos \theta_2  \\
		\end{pmatrix}
		\end{eqnarray} 
		In addition, we know that $\Sigma_{g_1 g_2} = \Sigma_{g_1 g_2}^{'}$, which implies that 
		\begin{eqnarray}
		\begin{pmatrix} 
		1 & 0  \\
		0 & \sigma_{g_1 g_2} \\
		\end{pmatrix}
		= 
		\begin{pmatrix} 
		\cos \theta_1 & \sin \theta_1  \\
		- \sin \theta_1 & \cos \theta_1  \\
		\end{pmatrix}
		\begin{pmatrix} 
		1 & 0  \\
		0 & \sigma_{g_1 g_2}^{'} \\
		\end{pmatrix}
		\begin{pmatrix} 
		\cos \theta_2 & - \sin \theta_2  \\
		\sin \theta_2 & \cos \theta_2  \\
		\end{pmatrix}.
		\end{eqnarray}
		After simplification, we have that 
		\begin{eqnarray}\label{tan}
		0 &=& -\cos \theta_1 \sin \theta_2 + \sigma^{'} \sin \theta_1 \cos \theta_2 \\
		0 &=& - \sin \theta_1 \cos \theta_2 + \sigma^{'} \cos \theta_1 \sin \theta_2.
		\end{eqnarray}
		Observe that $\cos \theta_1$ and $\cos \theta_2$ are not equal to zero, otherwise $\cos \theta_1 \cos \theta_2 + \sigma^{'} \sin \theta_1 \sin \theta_2 < 1$.
		By \eqref{tan}, we have $\tan \theta_1 = \sigma^{' 2}\tan \theta_1$. This implies that $\theta_1 = 0, ~ \textrm{or} ~ \pi$. It implies that 
		$\bar A_{g_1} = \bar A_{g_1}^{'}$ and $\sigma_{g_1 g_2} = \sigma_{g_1 g_2}^{'}$.
		
		Take any $g$ not in $\H{3}$ and $g_1$ in $\H{3}$, we know that $\bar A_g$ can be represented as $\mathbf a_g[\mathcal B_g] (c_1, c_2)$. It is easy to see that 
		$A_g^{'} = \mathbf a_g[\mathcal B_g] (c_1^{'}, c_2^{'})$ with $c_1^{'2} + c_2^{'2} = c_1^2 + c^2$. Compare $\Sigma_{g_1 g}$ and $\Sigma_{g_1 g}^{'}$, we have that 
		\begin{eqnarray}
		\begin{pmatrix} 
		1 & 0  \\
		0 & \sigma_{g_1 g_2} \\
		\end{pmatrix}
		\begin{pmatrix} 
		c_1 \\
		c_2 \\
		\end{pmatrix}
		= 
		\begin{pmatrix} 
		1 & 0  \\
		0 & \sigma_{g_1 g_2}^{'} \\
		\end{pmatrix}
		\begin{pmatrix} 
		c_1^{'} \\
		c_2^{'} \\
		\end{pmatrix}.
		\end{eqnarray}
		This implies that $c_1^{'} = c_1$, $c_2^{'} =  c_2$ and $\sigma_{g_1g}^{'} = \sigma_{g_1g}$.
		
		Furthermore, if $G > |\H{3}| $, we take any testlet pair $g_1$ and $g_2$ not in $\H{3}$. By $\Sigma_{g_1 g_2} = \Sigma_{g_1 g_2}$, we have 
		\begin{eqnarray}
		\begin{pmatrix} 
		c_{g_1 1} & c_{g_1 2} \\
		\end{pmatrix} 
		\begin{pmatrix} 
		1 & 0  \\
		0 & \sigma_{g_1 g_2} \\
		\end{pmatrix}
		\begin{pmatrix} 
		c_{g_2 1} \\
		c_{g_2 2} \\
		\end{pmatrix} 
		= 
		\begin{pmatrix} 
		c_{g_1 1} & c_{g_1 2} \\
		\end{pmatrix}
		\begin{pmatrix} 
		1 & 0  \\
		0 & \sigma_{g_1 g_2}^{'} \\
		\end{pmatrix}
		\begin{pmatrix} 
		c_{g_2 1} \\
		c_{g_2 2} \\
		\end{pmatrix}, 
		\end{eqnarray}
		which implies that $\sigma_{g_1 g_2}^{'} = \sigma_{g_1 g_2}$.
		Hence, all parameters are identifiable. This concludes our proof.
	\end{proof}
	
	\bigskip
	
	\begin{proof}[Proof of Theorem \ref{extendedH2}]
		For simplicity, we suppose $|\H{3}| = 1$ and $g_1 \in \H{3}$. For any $g \neq g_1$, it must hold that
		\begin{eqnarray}
		\begin{pmatrix} 
		1 & 0  \\
		0 & \sigma_{g g_1}^{'} \\
		\end{pmatrix}
		\begin{pmatrix} 
		\cos \theta_g & \sin \theta_g  \\
		- \sin \theta_g & \cos \theta_g \\
		\end{pmatrix}
		\begin{pmatrix} 
		c_{g1}  \\
		c_{g2} \\
		\end{pmatrix}
		= 
		\begin{pmatrix} 
		\cos \theta_1 & \sin \theta_1  \\
		- \sin \theta_1 & \cos \theta_1 \\
		\end{pmatrix}
		\begin{pmatrix} 
		1 & 0  \\
		0 & \sigma_{g g_1} \\
		\end{pmatrix}
		\begin{pmatrix} 
		c_{g1}  \\
		c_{g2} \\
		\end{pmatrix}
		\end{eqnarray} 
		according to Lemma \ref{fund:lem}, where $\bar A_g$ has the form of $ \mathbf a_g[\mathcal B_g] (c_{g1}, c_{g2})$. Here $c_2$ is a non-zero constant. 
		By simplification, we then have
		\begin{eqnarray}\label{implicit:eqn}
		\begin{pmatrix} 
		c_{g1} \cos \theta_g +  c_{g2} \sin \theta_g  \\
		\sigma_{g g_1}^{'} (- c_{g1} \sin \theta_g + c_{g2} \cos \theta_g) \\
		\end{pmatrix}
		=
		\begin{pmatrix} 
		c_{g1} \cos \theta_1 +  c_{g2} \sigma_{gg_1} \sin \theta_1  \\
		(- c_{g1} \sin \theta_1 + c_{g2} \sigma_{gg_1} \cos \theta_1) \\
		\end{pmatrix}.
		\end{eqnarray}
		Equation \eqref{implicit:eqn} can be viewed as a function of $\theta_1$, $\theta_g$ and $\sigma_{gg_1}^{'}$.
		Clearly, it admits the solution $\theta_g =0, \theta_1 = 0, \sigma_{gg_1}^{'} = \sigma_{gg_1}$. 
		
		We perturb $\theta_1$ a bit at $\theta_1 = 0$ locally, i.e., $\theta_1 = \delta$.
		Here we can always choose $\delta$ such that it keeps the sign of the first item in $\mathcal B_g$ to be positive. Otherwise, we can choose $-\delta$.
		Thus, by the implicit function theorem, \eqref{implicit:eqn} admits the solution $\theta_g^{'} = \theta_g(\delta)$ and $\sigma_{gg_1}^{'} = \sigma_{gg_1}(\delta)$, since the determinant of gradient does not vanish at $\theta_g = 0, \sigma_{gg_1}^{'} = \sigma_{gg_1}$ when $c_{g2} \neq 0$. 
		In addition, we know that $\theta_g(\delta) \rightarrow 0$ as $\delta \rightarrow 0$. Then it keeps the sign of first non-zero item in $\mathcal B_g$ since both $a_0[j_g]$ and $a_g[j_g] > 0$ when $\delta$ is close enough to $0$. (Here $j_g$ is the first non-zero item in $\mathcal B_g$.)
		Furthermore, it admits that $\sigma_{gg^{'}}^{'} = \sigma_{gg^{'}}(\delta)$ for any item pair $g, g^{'} \neq g_1$. 
		Let $\delta$ go to 0, then $\sigma_{g g^{'}}^{'}$s uniformly go to $\sigma_{gg^{'}}$ since the number of parameters is finite. By eigenvalue perturbation theory, there exists a $\delta$ such that that $\Sigma^{'}$ is still positive definite. 
		This guarantees that $\Sigma^{'}$ is still a covariance matrix. 
		
		By the above displays, the model is not identifiable since we have constructed another set of parameters leading to the same distribution.
		By the same technique, the model is not identifiable when $|\H{3}| = 0$. Thus we conclude the proof. 
	\end{proof}
	
	\bigskip
	
	Now we are ready to prove the main results.
	
	\begin{proof}[Proof of Theorem \ref{thm:extendunknown}]
		\textbf{Proof under Condition $\mathit{E2S}$}:
		The first step is to show that the covariance is identifiable.
		We prove this by checking the Condition C0 of Theorem \ref{thm:generalother}.
		It is also suffices to check the Condition C1 of Theorem \ref{thm:generalcov}.  
		
		By the requirement that $|\mathcal Q_g| \geq 3$, we denote these three items in testlet $g$ as $j_{g,1}, j_{g,2}$ and $j_{g,3}$. 
		From the requirement that $|\H{2}|\geq 1$, we can assume $g_1 \in \H{2}$ and $\bar A_{g_1}[\mathcal B_{g_1},:]$, $\bar A_{g_1}[\mathcal B_{g_2},:]$ are full-column rank with $\mathcal B_{g_1,1} = \{j_{g_1,1}, j_{g_1,2}\}$ and $\mathcal B_{g_1,2} = \{j_{g_1,3}, j_{g_1,4}\}$. 
		By the requirement of $|\H{3}|\geq 2$, we know there exists $g_2 \in \H{3}$ such that $g_2 \neq g_1$. 
		Then we can assume $\bar A_{g_2}[j_{g_2,2},:]$ and $\bar A_{g_2}[j_{g_2,3},:]$ are linearly independent for items $j_{g_2,2}$ and $j_{g_2,3}$. 
		
		We then can construct a partition $\mathcal B_1 \cup \mathcal B_2$ satisfying that $\mathcal B_1 = \{j_{g_1,1}, j_{g_1,2}, j_{g_2,1}, \ldots, j_{g,1};  g \neq g_1, g \neq g_2\}$, $\mathcal B_2 = \{j_{g_1,3}, j_{g_1,4}, j_{g_2,2}, j_{g_2,3}, \ldots, j_{g,2}, j_{g,3}; g \neq g_1, g \neq g_2\}$.
		It is easy to check that $A[\mathcal B_1,:]$ has full column rank. 
		Let $\mathcal B_{2a} = \{ j_{g_1,3}, j_{g_1,4}, j_{g_2,2}, \ldots, j_{g,2}; ; g \neq g_1, g \neq g_2 \}$. It is also easy to check that $A[\mathcal B_{2a},:]$ has full column rank and $A[\mathcal B_{2} - \{j\},:]$ has full column rank for $\forall j \in \mathcal B_{2a}$. Thus Condition C1 is satisfied.
		
		Hence the problem is reduced to the linear case with known variance.
		We only need to check the condition that $|\mathcal N| = 0$ according to Theorem \ref{extended:suff:general}.
		If not, there is a $g$ such that $\mathbf a_g = \mathbf 0$. It contradicts with $|\mathcal Q_g|\geq 3$. This concludes the proof.
		
		\textbf{Proof under Condition $\mathit{E1S}$}:
		The first step is still to show that the item covariance matrix is identifiable.
		Again, we prove this by checking the Condition C1. 
		%of Theorem \ref{thm:generalother}.
		%It is also suffices to check the Condition C1 of Theorem \ref{thm:generalcov}.  
		
		By the requirement that $|\mathcal Q_g| \geq 3$ for each testlet $g$, we take any three items in testlet $g$ and denote them as $j_{g,1}, j_{g,2}$ and $j_{g,3}$. 
		By the requirement of $|\H{3}|\geq 3$, we can assume $g_1, g_2, g_3 \in \H{3}$ and assume $\bar A_{g_1}[\mathcal B_{g_1,a},:]$, $\bar A_{g_2}[\mathcal B_{g_2,a},:]$, $\bar A_{g_3}[\mathcal B_{g_3,a},:]$ have full-column rank with $\mathcal B_{g_1,a} = \{j_{g_1,1}, j_{g_1,2}\}$, $\mathcal B_{g_2,a} = \{j_{g_2,1}, j_{g_2,2}\}$ and $\mathcal B_{g_3,a} = \{j_{g_3,1}, j_{g_3,2}\}$. 
		
		In the following, we need to verify that $\{1,\ldots,J\} - j$ can be partitioned into two item sets $\mathcal B_{1,j}$ and $\mathcal B_{2,j}$ for each item $j$ such that $A[\mathcal B_{1,j},:]$ and $A[\mathcal B_{2,j},:]$ satisfy Condition C1. 
		
		If $j$ belongs to testlet $g$ ($g \neq g_1, g_2, g_3$), we can assume $j = j_{g,1}$ without loss of generality. Then we can set 
		$$\mathcal B_{1,j} = \{j_{g_1,1}, j_{g_1,2}, j_{g_2,3}, j_{g_3,3}, \ldots, j_{g,2}, \ldots;  g \neq g_1, g_2, g_3 \},$$
		and 
		$$\mathcal B_{2,j} = \{j_{g_1,3}, j_{g_2,1}, j_{g_2,2}, j_{g_3,1}, j_{g_3,2}, \ldots, j_{g,3}, \ldots; g \neq g_1, g_2, g_3 \}.$$
		
		If $j$ belongs to testlet $g$ ($g \in \{g_1, g_2, g_3\}$), we can assume $j = j_{g_1,1}$ without loss of generality. Then we can set 
		$$\mathcal B_{1,j} = \{j_{g_1,2}, j_{g_2,1}, j_{g_2,2}, j_{g_3,3}, \ldots, j_{g,2}, \ldots;  g \neq g_1, g_2, g_3 \},$$
		and 
		$$\mathcal B_{2,j} = \{j_{g_1,3}, j_{g_2,3}, j_{g_3,1}, j_{g_3,2}, \ldots, j_{g,3}, \ldots; g \neq g_1, g_2, g_3 \}.$$
		Then by Condition C1 we know that the covariance matrix is identifiable.
		Obviously $|\mathcal N| = 0$ and $|\H{3}| \geq 2$ hold. Hence we conclude the proof.
	\end{proof}
	
	\bigskip
	
	\begin{proof}[Proof of Theorem \ref{thm:extendunknown:nec}]
		To prove the necessity that $|\H{3}| \geq 2$, we want to show that we can always construct another set of parameters leading to the same distribution for any $\Theta$ satisfies $\H{1} \leq 1$. 
		
		Under this case, we keep $\boldsymbol \lambda$ fixed and only consider to construct another loading matrix $A$. 
		It is easy to compute that the covariance between items $j_1$ and $j_2$ from group $g$ is
		\begin{eqnarray}
		\sigma_{j_1 j_2} = \mathbf a_0 [j_1] \mathbf a_0 [j_2] + \mathbf a_g [j_1] \mathbf a_g [j_2];
		\end{eqnarray}
		the covariance between items $j_1$ and $j_2$ from groups $g_1$ and $g_2$ is
		\begin{eqnarray}
		\sigma_{j_1 j_2} = \mathbf a_0 [j_1] \mathbf a_0 [j_2] + \mathbf a_{g_1} [j_1] \mathbf a_{g_2} [j_2] \sigma_{g_1 g_2}.
		\end{eqnarray} 
		We can write them in matrix form which becomes
		\begin{eqnarray}
		\Sigma_{gg} = \bar A_g \bar A_g^T ~ \textrm{and} ~ \Sigma_{g_1 g_2} = \bar A_{g_1} \begin{pmatrix} 
		1 & 0 \\
		0 & \sigma_{g_1 g_2} \\
		\end{pmatrix}
		\bar A_{g_2}^T.
		\end{eqnarray}
		Then, by Theorem \ref{extendedH2}, we have that $|\H{3}| \geq 2$.
		
		To prove the necessity that $|\mathcal Q_g| \geq 2$ for all $g \in \{1, \ldots, G\}$, we construct another set of parameters leading to the same distribution.
		Take any $g$ such that $|\mathcal Q_g| \leq 1$. Without loss of generality, we can assume $|\mathcal Q_g| = 1$ and assume item $j$ satisfies $\mathbf a_g[j] \neq 0$. Thus, we can construct another set of parameters such that $\mathbf a_g^{'}[j] = x \cdot a_g[j]$, $\Sigma_G^{'}[g,g^{'}] = \Sigma_G[g,g^{'}]/x$ for $g^{'} \neq g$ and keep other parameters fixed. It can be verified that this set of parameters works. 
		
		Furthermore, by the same argument in the proof of Theorem \ref{thm:unknown:bifactor}, we must have $|\mathcal Q_g| \geq 3$ for all $g$ satisfying $\Sigma_G[g, -g] = \mathbf 0$,
	\end{proof}

	\bigskip
	
	\begin{proof}[Proof of Proposition \ref{prop:extend:H5:unnec}]
		We still use proof by contradiction.
		The first step is to show that the main factor loadings are identifiable. 
		By assumption, we must have that 
		\begin{eqnarray}\label{H5:identity}
		\bar A_{g_1}^{'} \begin{pmatrix}
		1 & 0 \\
		0 & \Sigma_G^{'}[g_1,g_2] \\
		\end{pmatrix}
		\bar A_{g_2}^{'T}
		= 
		\bar A_{g_1} \begin{pmatrix}
		1 & 0 \\
		0 & \Sigma_G[g_1,g_2] \\
		\end{pmatrix}
		\bar A_{g_2}^T.
		\end{eqnarray}
		Since $g_1, g_2 \in \H{3}$, we know $\bar A_{g_1}, \bar A_{g_2}$ are full column rank matrices. There must exist full rank matrices $R_1$ and $R_2$ such that $\bar A_{g_1}^{'} = \bar A_{g_1} \cdot R_1$ and $\bar A_{g_2}^{'} = \bar A_{g_2} \cdot R_2$.
		
		Next, again by assumption, we know that 
		\begin{eqnarray}
		\bar A_{g_1}^{'}[-j,:] \bar A_{g_1}^{'}[j,:]^T
		= 
		\bar A_{g_1}[-j,:] \bar A_{g_1}[j,:]^T; ~~ j \in \mathcal Q_{g_1}.
		\end{eqnarray}
		This implies that $R_1 R_1^T A_{g_1}[j,:]^T = A_{g_1}[j,:]^T$ for $j \in \mathcal Q_{g_1}$, as the Kruskal rank of $\bar A_{g_1}$ is 2. It must hold that $R_1 R_1^T = I$. We then parameterize $R_1$ as
		\begin{eqnarray}
		R_1 = \begin{pmatrix}
		\cos \theta & - \sin \theta \\
		\sin \theta  & \cos \theta \\
		\end{pmatrix}.
		\end{eqnarray}
		By \eqref{H5:identity}, we further have 
		\begin{eqnarray}
		R_2 = 
		\begin{pmatrix}
		1 & 0\\
		0 & \frac{1}{\sigma^{'}}\\
		\end{pmatrix}
		\begin{pmatrix}
		\cos \theta & \sin \theta \\
		-\sin \theta & \cos \theta\\
		\end{pmatrix}
		\begin{pmatrix}
		1 & 0\\
		0 & \sigma\\
		\end{pmatrix}.
		\end{eqnarray}
		Again by assumption, we know that 
		\begin{eqnarray}
		\bar A_{g_2}^{'}[-j,:] \bar A_{g_2}^{'}[j,:]^T
		= 
		\bar A_{g_2}[-j,:] \bar A_{g_2}[j,:]^T; ~~ j \in \mathcal Q_{g_2}.
		\end{eqnarray}
		This implies that determinant of $R_2 R_2^T - I$ is zero. By straightfoward calculation,
		\begin{eqnarray*}
			0 & = & \mathrm{det}(R_2 R_2^T - I) \\
			& = & \mathrm{det}(
			\begin{pmatrix}
				\sin^2 \theta \cdot (\frac{1}{\sigma^{' 2}} - 1) & \sigma \cos \theta \sin \theta \cdot (1 - \frac{1}{\sigma^{' 2}})\\
				\sigma \cos \theta \sin \theta \cdot (1 - \frac{1}{\sigma^{' 2}}) & \frac{\sigma^2}{\sigma^{'2}} \cos^2 \theta + \sigma^2 \sin^2 \theta - 1\\
			\end{pmatrix}
			) \\
			& = & \sin^2 \theta \cdot (\frac{1}{\sigma^{' 2}} - 1) \cdot (\sigma^2 - 1).
		\end{eqnarray*}
		Hence, $\theta = 0$ and $R_1 = I$. This implies that main factor loading $\mathbf a_{g_1}$ is identifiable.
		By assumption, it holds that 
		\begin{eqnarray}
		\bar A_{g_1}
		\begin{pmatrix}
		1 & 0 \\
		0 & \Sigma_G[g_1,g]^{'} \\
		\end{pmatrix}
		\bar A_{g}^{'T}
		= 
		\bar A_{g_1}
		\begin{pmatrix}
		1 & 0 \\
		0 & \Sigma_G[g_1,g]^{} \\
		\end{pmatrix}
		\bar A_{g}^{T}
		\end{eqnarray}
		for any $g \neq g_1$.
		We must have that $\mathbf a_g^{'} = \mathbf a_g^{}$, since $\bar A_{g_1}$ has full column rank. Hence all main factor loadings are identifiable. 
		
		Next, we prove the identifiability of testlet-specific loadings.
		By assumption we have that 
		\begin{eqnarray}\label{H5:Q2}
		\mathbf a_g^{'}[-j] \cdot \mathbf a_g^{'}[j]
		= \mathbf a_g^{}[-j] \cdot \mathbf a_g^{}[j]
		\end{eqnarray}
		for any $j$ and $g$. If $|\mathcal Q_g| \geq 3$, we have that $\mathbf a_g = \mathbf a_g^{'}$ according to proof in Theorem \ref{thm:unknown:bifactor}.
		If $|\mathcal Q_g| = 2$, we know that $\Sigma_G[g,-g]$ is not zero-vector. Hence, we must have that $\mathbf a_g^{'} = x \cdot \mathbf a_g$ for some positive $x$ by comparing the cross variance. Finally, by \eqref{H5:Q2}, we have $x^2 = 1$. Thus $x \equiv 1$. Then $\mathbf a_g$ is identifiable. We conclude the proof. 
	\end{proof}

	To end this section, we provide a variant of Theorem 5.1 in \cite{anderson1956statistical} to item factor models with probit link, which provides a basis to the identifiability results for the probit extended bifactor and two-tier models. Based on this result, there will be no difference between the proofs of linear bifactor models and probit bifactor models.
	\begin{theorem}\label{thm:generalcov}
		Under the general factor model with probit link as in \eqref{probit}, assume the following holds:
		\begin{itemize}
			\item[C1] There exists a partition of the items $\{1,\ldots,J\} = \mathcal B_1 \cup \mathcal B_2$ such that (1) $A[\mathcal B_1,:]$ is full-column rank; (2) there exists a subset of $\mathcal B_2$, $\mathcal B_{2a}$, satisfying that $A[\mathcal B_{2a},:]$ is full-column rank and  $A[\mathcal B_{2} - \{j\},:]$ is full-column rank for $\forall j \in \mathcal B_{2a}$. 
			%\item[C2'] $|Q_g| \geq 2$ for all $g$. For any $g \notin \H{3} \cap \mathcal H_7$, there exists $g_1 \in \mathcal \H{3} \cap \mathcal H_7$ such that $\sigma_{g g_1} \neq 0$.
		\end{itemize}
		Then, the covariance matrix $A \Sigma A^T$ is identifiable.
	\end{theorem}
	The result is closely related to Kruskal rank. That is, Condition C1 is satisfied when
	$\mathcal B_2$ contains a subset of $\mathcal B_{2a}$ of $K+1$ items and $A[\mathcal B_{2a},:]^T$ has Kruskal rank $K$. Thus, Theorem \ref{thm:generalcov} requires $2K + 1$ or more items.
	The condition is very weak, especially in terms of number of required items.
	Notice that there are $J K$ loading parameters and $J(J-1)/2$ restrictions. Then $J(J-1)/2 \geq J K$ if and only if $J \geq 2 K + 1$. 
	In fact, we can prove that $2 K + 1$ is minimal possible number of items for model identifiability for $K = 1$ and $2$.
	
	For the probit bifactor model, due to its special structure, the conditions for identifiability can be less restrictive than what Theorem \ref{thm:generalcov} requires for general probit factor models: Specifically, Theorem \ref{thm:generalcov} requires at least $3K + 1$ items for meeting Condition C1, which is slightly stronger than the optimal conditions for the bifactor model as stated in Theorem \ref{thm:probit:bifactor}. 
	
	\bigskip
	
	\begin{proof}[Proof of Theorem \ref{thm:generalcov}]
		By comparing the correlation between $\mathcal B_1$ and $\mathcal B_2$, we have that 
		\begin{eqnarray}\label{general:compare}
		\tilde A[\mathcal B_1,:] \Sigma \tilde A[\mathcal B_2,:]^T = \tilde A^{'}[\mathcal B_1,:] \Sigma^{'} \tilde A^{'}[\mathcal B_2,:]^T.
		\end{eqnarray} 
		Since both $\tilde A[\mathcal B_1,:]$ and $\tilde A[\mathcal B_2,:]$ have full column rank, we can assume $\tilde A^{'}[\mathcal B_1,:] = \tilde A[\mathcal B_1,:] R_1$ and $\tilde A^{'}[\mathcal B_2,:] = \tilde A[\mathcal B_2,:] R_2$. Thus \eqref{general:compare} becomes
		\begin{eqnarray}\label{general:compare:reduced}
		\Sigma  = R_1 \Sigma^{'} R_2^T.
		\end{eqnarray}
		
		By comparing the correlation within $\mathcal B_2$, we have that 
		\begin{eqnarray}\label{general:compare:within}
		\tilde A[\mathcal B_2 - \{j\},:] \Sigma \tilde A[j,:]^T = \tilde A^{'}[\mathcal B_2 - \{j\},:] \Sigma^{'} \tilde A^{'}[j,:]^T \quad \forall j \in \mathcal B_{2a}.
		\end{eqnarray}
		Notice that $\tilde A[\mathcal B_2 - \{j\},:]$ has full rank $K$ as $A[\mathcal B_2 - \{j\},:]$ does. \eqref{general:compare:within} is reduced to 
		\begin{eqnarray}
		\Sigma \tilde A[j,:]^T = R_2 \Sigma^{'} R_2^T \tilde A[j,:]^T \quad \forall j \in \mathcal B_{2a},
		\end{eqnarray}
		that is $ (\Sigma - R_2 \Sigma^{'} R_2^T)\tilde A[j,:]^T = \mathbf 0$ for all $j \in \mathcal B_{2a}$. This implies that the kernel of the linear map $\Sigma - R_2 \Sigma^{'} R_2^T$ is the whole vector space. Thus $\Sigma \equiv R_2 \Sigma^{'} R_2^T$ which further implies $R_1 = R_2$. Therefore,
		\begin{eqnarray}
		\tilde A \Sigma \tilde A^T = \tilde A^{'} \Sigma^{'} \tilde (A^{'})^T.
		\end{eqnarray} 
		We conclude that the covariance matrix is identifiable by comparing the diagonals of the above equation.
	\end{proof}

	\newpage
	
	\section{Proofs for Two-Tier Models}\label{sec:E}
	
	Following similar rationale as for the extended bifactor models, it suffices to prove the following three theorems.
	
	\begin{theorem}\label{twotier:T1}
		Under a linear two-tier model, suppose $\boldsymbol \lambda$ is known. Then the parameter is identifiable if it satisfies $\mathit{T1S}$.
	\end{theorem}
	
	\begin{theorem}\label{twotier:T2}
		Under a linear two-tier model, suppose $\boldsymbol \lambda$ is known. Then the parameter is identifiable if it satisfies $\mathit{T2S}$.
	\end{theorem}
	
	\begin{theorem}\label{twotier:T3}
		Under a linear two-tier model, suppose $\boldsymbol \lambda$ is known. Then the parameter is identifiable if it satisfies $\mathit{T3S}$.
	\end{theorem}

	\begin{proof}[Proof of Theorem \ref{twotier:T1}]
		We prove this by using method of contradiction. 
		%		Let $A_g$ be the row block of main factors corresponding to Testlet $g$ and $\mathbf a_g$ be the testlet factor corresponding to Testlet $g$. 
		%		Suppose there exists another set of parameter leading to the same distribution. 
		We adopt notation $\Sigma_{\mathcal H_1, \mathcal H_2}$ to denote the item covariance across testlets from $\mathcal H_1$ and $\mathcal H_2$.
		Since $|\mathcal H_3| \geq 3$, we then take three Testlets $g_1$, $g_2$ and $g_3$ such that $g_1, g_2, g_3 \in \mathcal H_3$. Note that $\Sigma_{g g^{'}}^{'} = \Sigma_{g g^{'}}$ by comparing the item covariance across different testlets.
		We then know that $A^{'}[\mathcal B_{g},1:L] = A_g [\mathcal B_{g},1:L] R_g$ ($g = g_1, g_2, g_3$), where $R_g$ is a $L$ by $L$ full rank matrix. We also know that $R_g^{'} \Sigma_L^{'} R_g = \Sigma_L$. This indicates that $R_{g_1} = R_{g_2} = R_{g_3} = R$ and $R \Sigma_L^{'} R^T =\Sigma_L$.
		It further implies 
		\begin{eqnarray}
		A[\mathcal B_{\mathcal H_3},1:L] R \Sigma_L^{'} R^T A[\mathcal B_{\mathcal H_3},1:L]^T = A[\mathcal B_{\mathcal H_3},1:L] \Sigma_L A[\mathcal B_{\mathcal H_3},1:L]^T
		\end{eqnarray}
		Since $A[\mathcal B_{\mathcal H_3},1:L]$ contains an identity, we then have $R = I$ and $\Sigma_{L}^{'} = \Sigma_L$. 
		
		For any $g \notin \mathcal H_3$, we know that 
		$\Sigma_{\mathcal H_3, g}^{'} = \Sigma_{\mathcal H_3, g}$, indicating that $A^{'}[\mathcal B_{g},1:L] = A[\mathcal B_{g},1:L]$. Lastly, we then have $\mathbf a_g^{'} = \mathbf a_g$ for all $g$ by comparing the item covariance within testlets.    	
	\end{proof}
	
	\bigskip
	
	\begin{proof}[Proof of Theorem \ref{twotier:T2}]
		We still use method of contradiction to prove this theorem. 	
		Suppose there exists another set of parameter leading to the same distribution. Since $|\mathcal H_3| \geq 2$ and $|\mathcal H_4| \geq 1$, we then take two Testlets $g_1$ and $g_2$ such that  $g_1 \in \mathcal H_3 \cap \mathcal H_4$ and $g_2 \in \mathcal H_3$. Notice that $\Sigma_{g_1 g_2}^{'} = \Sigma_{g_1 g_2}$.
		We then know that $A^{'}[\mathcal B_{g},1:L] = A[\mathcal B_{g},1:L] R_g$ ($g = g_1, g_2$), where $R_g$ is a $L$ by $L$ full rank matrix. 
		%	Hence, $R_{k_1} = R_{k_2} = R_{k_3}$, which implies that $R_k \Sigma_L^{'} R_k^T = \Sigma_L$. 
		Recall $\bar A_{g_1} = (A[\mathcal B_{g_1},1:L], \mathbf a_{g_1}[\mathcal B_{g_1}])$, we then know $\bar A_{g_1}^{'} = \bar A_{g_1} \begin{pmatrix}
		R_{g_1} & c \\
		0 & d \\
		\end{pmatrix}$.
		Furthermore, we know that $\Sigma_{g_1 g_1}^{'} = \Sigma_{g_1 g_1}$. This gives us that 
		\begin{eqnarray}
		c \cdot d = 0 ~ \textrm{and} ~ d^2 = 1, 
		\end{eqnarray}
		which implies that $\mathbf a_{g_1}^{'} (\mathbf a_{g_1}^{'})^T = \mathbf a_{g_1} \mathbf a_{g_1}^T$. 
		We then also have $R_{g_1} \Sigma_L^{'} R_{g_1}^T = \Sigma_L$.
		Since $R_{g_1} \Sigma_L^{'} R_{g_2}^T = \Sigma_L$, it implies that $R_{g_1} = R_{g_2}$. Then $R_g \Sigma_L^{'} R_g^T = \Sigma_L$ for any $g \in \mathcal H_3$.
		%Further notice that $\Sigma_{\mathcal H_3, \mathcal H_3}^{'} = \Sigma_{\mathcal H_3, \mathcal H_3}$, 
		It implies that 
		\begin{eqnarray}
		A[\mathcal B_{\mathcal H_3},1:L] R \Sigma_L^{'} R^T A[\mathcal B_{\mathcal H_3},1:L]^T = A[\mathcal B_{\mathcal H_3},1:L] \Sigma_L A[\mathcal B_{\mathcal H_3},1:L]^T
		\end{eqnarray}
		Since $A[\mathcal B_{\mathcal H_3},1:L]$ contains an identity, we then have $R = I$ and $\Sigma_{L}^{'} = \Sigma$. 
		
		Furthermore, for any $g \notin \mathcal H_3$, we know that 
		$\Sigma_{\mathcal H_3, g}^{'} = \Sigma_{\mathcal H_3, g}$, indicating that $A^{'}[\mathcal B_{g},1:L] = A[\mathcal B_{g},1:L]$. Lastly, we have $\mathbf a_g^{'} = \pm \mathbf a_g$ by using the fact that 
		$\Sigma_{gg}^{'} = \Sigma_{gg}$ for all $g$.    
	\end{proof}
	
	\bigskip
	
	\begin{proof}[Proof of Theorem \ref{twotier:T3}]
		By the same strategy, 
		suppose the model is not identifiable. There exists another set of parameters $A^{'}$ and $\Sigma^{'}$ leading to the same model. 
		Hence, $A^{'}[\mathcal B_{\mathcal G_1}, 1:L] \Sigma_L^{'} A^{'}[\mathcal B_{\mathcal G_2}, 1:L]^T = A^{}[\mathcal B_{\mathcal G_1}, 1:L] \Sigma_L A^{}[\mathcal B_{\mathcal G_2}, 1:L]^T$. Since both 
		$A[\mathcal B_{\mathcal G_1}, 1:L]$ and $A^{}[\mathcal B_{\mathcal G_2}, 1:L]$ are full column rank matrices according to Condition $\mathit{T3S}$-(b). This implies that $A^{'}[\mathcal B_{\mathcal G_1}, 1:L]$ spans the same subspace as $A[\mathcal B_{\mathcal G_1}, 1:L]$ does. In other words, $A^{'}[\mathcal B_{\mathcal G_1}, 1:L]$ has the form of $A[\mathcal B_{\mathcal G_1}, 1:L] \cdot D$ for some $D$ being a $L$ by $L$ full rank matrix.
		%	Further notice $I(A^{'}_{I_{\mathcal K_1}, 1:L}) = I(A_{I_{\mathcal K_1}, 1:L})$, this implies that $R$ should be a diagonal matrix.
		
		By Condition $\mathit{T3S}$-(a) that $span(A[\mathcal B_{\mathcal G_1}, 1:L]) \cap span(A^{'}[\mathcal B_{\mathcal G_1}, L + \mathcal G_1 ]) = \mathbf 0$ and 
		by assumption that $A^{'}[\mathcal B_{\mathcal G_1}, :] \Sigma^{'} A^{'}[\mathcal B_{\mathcal G_1}, :]^T = A^{}[\mathcal B_{\mathcal G_1}, :] \Sigma^{} A^{}[\mathcal B_{\mathcal G_1}, :]^T$, we then know $A^{'}[\mathcal B_{\mathcal G_1}, (1:L) + \mathcal G_1] = A^{}[\mathcal B_{\mathcal G_1}, (1:L) + \mathcal G_1] \cdot R$ where $R$ is $L+G_1$ by $L+G_1$ matrix (Here $G_1 = |\mathcal G_1|$). Furthermore, $R$ should has the form of $\begin{pmatrix}
		D & X_1 \\
		\mathbf 0 & X_2 \\
		\end{pmatrix}$.
		As a result, we then have 
		\begin{eqnarray}
		\begin{pmatrix}
		D & X_1 \\
		\mathbf 0 & X_2 \\
		\end{pmatrix}
		\begin{pmatrix}
		\Sigma_L^{'} & 0 \\
		\mathbf 0 & I \\
		\end{pmatrix}
		\begin{pmatrix}
		D & \mathbf 0 \\
		X_1^T & X_2^T \\
		\end{pmatrix}
		= 
		\begin{pmatrix}
		\Sigma_L & \mathbf 0 \\
		\mathbf 0 & I \\
		\end{pmatrix},	
		\end{eqnarray}
		which implies that $X_1 X_2^T = \mathbf 0$, $X_2 X_2^T = I$, $D \Sigma_L^{'} D + X_1 X_1^T = \Sigma_L$.
		Therefore, it gives $X_1 = \mathbf 0$, $X_2 = I$ and $D \Sigma_L^{'} D^T = \Sigma_L$.
		This further gives us that $A[:,1:L] D \Sigma_L^{'} D^T A[:,1:L]^T = A[:,1:L] \Sigma_L A[:,1:L]^T$. Since $A_{:,1:L}$ contains an identity, we then have $D$ is an identity matrix.
		Thus, we have $A[:,1:L]^{'} = A[:,1:L]$ which further implies $A^{'} = A$. We thus conclude the proof. 
	\end{proof}

	\newpage
	
	\section{Proof of Theorem \ref{bifactor:general:nonid} in Section 4}\label{sec:F}
	
	\begin{proof}[Proof of Theorem \ref{bifactor:general:nonid}]
		Let $A_G$ be $(\mathbf a_1, \ldots, \mathbf a_G)$. We prove the results by considering the following two cases.
		\begin{itemize}
			\item[] \textbf{Case 1}: $\mathbf a_0$ is not in the range of $A_G$.
			Under this case, we have that $\mathbf a_0, \ldots, \mathbf a_G$ are linearly independent. If not, there exists $c_0, \ldots, c_G$ such that $ \sum_{g=0}^G c_g \mathbf a_g = \mathbf 0$. Then we have  
			$\sum_{g=1}^G c_g \mathbf a_g = \mathbf 0$ by assumption that $a_0$ is not in the range of $A_0$.
			Since $\mathbf a_g$ has non-zero loadings on $g$th testlet, $\mathbf a_1, \ldots, \mathbf a_G$ are linearly independent. Thus we have $c_0 = \cdots = c_G =0$.
			
			If there is another model $\mathcal P^{'}(A^{'},d,\Sigma_G^{'},\rho^{'})$ that implies the same distribution, we must have
			\begin{eqnarray}
			A^{'} = A 
			\begin{pmatrix}
			x & \mathbf 0^T\\
			\mathbf y & \Lambda \\
			\end{pmatrix}
			\end{eqnarray} 
			since $A$ has full column rank. 
			Furthermore, we must have 
			\begin{eqnarray}
			\begin{pmatrix}
			x & \mathbf 0^T\\
			\mathbf y & \Lambda \\
			\end{pmatrix}
			\begin{pmatrix}
			1 & (\rrho^{'})^T \\
			\rrho^{'}  & \Sigma_G^{'} \\
			\end{pmatrix}
			\begin{pmatrix}
			x & \mathbf y^T\\
			\mathbf 0 & \Lambda \\
			\end{pmatrix}
			= 
			\begin{pmatrix}
			1 & \rrho^T \\
			\rrho  & \Sigma_G \\
			\end{pmatrix}.
			\end{eqnarray}
			This gives us that $x^2 = 1$, $x(\mathbf y + \Lambda \rrho^{'}) = \rrho$ and $\mathbf y \mathbf y^T + \Lambda \rrho^{'} \mathbf y^T + \mathbf y (\rrho^{'})^T \Lambda + \Lambda \Sigma_G^{'} \Lambda = \Sigma_G$.
			
			After simplification, we have 
			\begin{eqnarray*}
				(\mathbf y + \Lambda \rrho^{'}) (\mathbf y + \Lambda \rrho^{'})^T - (\Lambda \rrho^{'})(\Lambda \rrho^{'})^T
				+ \Lambda \Sigma_G^{'} \Lambda &=& \Sigma_G \\
				\rrho \rrho^T - (\Lambda \rrho^{'})(\Lambda \rrho^{'})^T
				+ \Lambda \Sigma_G^{'} \Lambda = \Sigma_G \\
				\rrho \rrho^T - (\tilde \rrho)( \tilde \rrho)^T
				+ \Lambda \Sigma_G^{'} \Lambda = \Sigma_G,
			\end{eqnarray*}
			where $\tilde \rrho = \Lambda \rrho^{'}$.
			
			We choose $\tilde \rrho$ to be arbitrary close to $\rrho$, set $\Lambda[g,g] = \sqrt{(\Sigma_G + (\tilde \rrho) (\tilde \rrho)^T - \rrho \rrho^T)[g,g]}$ and set $\Sigma_G^{'} = \Lambda^{-1}(\Sigma_G + (\tilde \rrho) (\tilde \rrho)^T - \rrho \rrho^T ) \Lambda^{-1}$. 
			Hence we find another set of parameters which leads to the same distribution. Thus the model is not identifiable.
			
			%		Set $x = 1$, $\mathbf y = 0$, $\Lambda \rrho^{'} = \rrho$, we then have
			%		$\Lambda \Sigma_G^{'} \Lambda = \Sigma_G$. 
			\item[] \textbf{Case 2}: $\mathbf a_0$ is in the range of $A_G$.
			We construct another model $\mathcal P(A^{'}, d^{'}, \Sigma_G^{'}, \rrho^{'})$ such that $A_G^{'} = A_G$, $\Sigma_G^{'} = \Sigma_G$.
			To determine the value of $\mathbf a_0^{'}$ and $\rrho^{'}$, we use the following equations.
			\begin{eqnarray*}
				(\mathbf a_0, A_G) 
				\begin{pmatrix}
					1 & \rrho \\
					\rrho & \Sigma_G \\
				\end{pmatrix}
				(\mathbf a_0, A_G)^T
				&=&
				(\mathbf a_0^{'}, A_G) 
				\begin{pmatrix}
					1 & \rrho^{'} \\
					\rrho^{'} & \Sigma_G \\
				\end{pmatrix}
				(\mathbf a_0^{'}, A_G)^T  \\
				\mathbf a_0 \mathbf a_0^T
				+ A_G \rrho \mathbf a_0^T + \mathbf a_0 \rrho^T A_G^T
				&=&
				(\mathbf a_0^{'}) (\mathbf a_0^{'})^T
				+ A_G \rrho^{'} (\mathbf a_0^{'})^T + (\mathbf a_0^{'}) (\rrho^{'})^T A_G^T \\
				(\mathbf a_0 + A_G \rrho)(\mathbf a_0 + A_G \rrho)^T - (A_G \rrho)(A_G \rrho)^T
				&=&
				(\mathbf a_0^{'} + A_G \rrho^{'})(\mathbf a_0^{'} + A_G \rrho^{'})^T - (A_G \rrho^{'})(A_G \rrho^{'})^T \\
				\mathbf x \mathbf x^T - \mathbf y \mathbf y^T
				&=& 
				(\mathbf x^{'}) (\mathbf x^{'})^T - (\mathbf y^{'}) (\mathbf y^{'})^T
			\end{eqnarray*} 
			where $\mathbf x = \mathbf a_0 + A_G \rrho$, $\mathbf y = A_G \rrho$ and $\mathbf x^{'}$, $\mathbf y^{'}$ are defined correspondingly.
			It is easy to check that $\mathbf x^{'} = \sqrt{1 + c^2} \mathbf x + c \mathbf y$, $\mathbf y^{'} = c \mathbf x + \sqrt{1 + c^2} \mathbf y$ satisfies the above equations. 	   
			This give us 
			\begin{eqnarray}\label{eqn:a}
			\mathbf a_0^{'} + A_G \mathbf \rrho^{'} = \sqrt{1 + c^2} (\mathbf a_0 + A_G \mathbf \rrho) + c A_G \rrho
			\end{eqnarray}
			and
			\begin{eqnarray}\label{eqn:rho}
			A_G \mathbf \rrho^{'} = c (\mathbf a_0 + A_G \mathbf \rrho) + \sqrt{1 + c^2} A_G \rrho
			\end{eqnarray} 
			By assumption that $\mathbf a_0$ is in the range of $A_G$, there exists $\mathbf b_0$ such that $\mathbf a_0 = A_G \mathbf b_0$. 
			Equation \eqref{eqn:rho} becomes $A_G \rrho^{'} = A_G (c \mathbf b_0 + (c + \sqrt{1 + c^2}) \rrho)$. This implies $\rrho^{'} = c \mathbf b_0 + (c + \sqrt{1 + c^2}) \rrho$. Plug this into \eqref{eqn:a}, we have 
			\begin{eqnarray*}
				\mathbf a_0^{'} + A_G (c \mathbf b_0 + (c + \sqrt{1 + c^2} \rrho)) = \sqrt{1 + c^2} (\mathbf a_0 + A_G \mathbf \rrho) + c A_G \rrho.
			\end{eqnarray*}
			Thus $\mathbf a_0^{'} = \sqrt{1 + c^2} (\mathbf a_0 + A_G \mathbf \rrho) + c A_G \rrho - A_G (c \mathbf b_0 + (c + \sqrt{1 + c^2}) \rrho) = (\sqrt{1 + c^2} - c)\mathbf a_0$.
			
			In summary, we constructed a different set of parameters which has the same distribution as the true model.
		\end{itemize}
		Hence we conclude the proof.
	\end{proof}

\bibliography{references.bib}

\end{document}